\documentclass{amsart}
\usepackage{amssymb, amsmath, tikz-cd, amsthm,amsfonts}
\usepackage{mathtools}

\newcommand{\nn}{\mathbb N}
 
\newcommand{\rr}{\mathbb R}

\newcommand{\bb}{\mathfrak{B}}

\newcommand{\norm}[1]{\left\lVert#1\right\rVert}

\newtheorem{thm}{Theorem}
\newtheorem{prop}{Proposition}
\newtheorem{lem}{Lemma}
\newtheorem{coro}{Corollary}

\newtheorem{exam}{Example}
\newtheorem{rem}{Remark}

\newtheorem{definition}{Definition}

\begin{document}

\title{ Induced Actions of $\bb$-Volterra Operators on Regular Bounded Martingale Spaces}
\author[Nazife Erkur\c{s}un-\"Ozcan, Niyazi An{\i}l Gezer ]{Nazife Erkur\c{s}un-\"Ozcan$^{(1)}$, Niyazi An{\i}l Gezer$^{(2, *)}$}

\address{$^{1}$ Department of Mathematics, Faculty of Science, Hacetttepe University, Ankara, 06800,Turkey.}
\email{{erkursun.ozcan@hacettepe.edu.tr}}

\address{$^{2}$ Independent Researcher, Ankara, Turkey.}
\email{{anilgezer@gmail.com}}



\date{Received: xxxxxx; Revised: yyyyyy; Accepted: zzzzzz.
\newline \indent $^{*}$ Corresponding author}

\begin{abstract} A positive operator $T:E\to E$ on a Banach lattice $E$ with an order continuous norm is said to be $\bb$-Volterra with respect to a Boolean algebra $\bb$ of order projections of $E$ if the bands canonically corresponding to elements of $\bb$ are left fixed by $T$. A linearly ordered sequence $\xi$ in $\bb$ connecting $\bold{0}$ to $\bold{1}$ is called a forward filtration. A forward filtration can be to used to lift the action of the $\bb$-Volterra operator $T$ from the underlying Banach lattice $E$ to an action of a new norm continuous operator $\hat{T}_{\xi}\colon \mathcal{M}_{r}(\xi) \to \mathcal{M}_{r}(\xi)$ on the Banach lattice $\mathcal{M}_{r}(\xi)$ of regular bounded martingales on $E$ corresponding to $\xi$. In the present paper, we study properties of these actions. The set of forward filtrations are left fixed by a function which erases the first order projection of a forward filtration and which shifts the remaining order projections towards $\textbf{0}$. This function canonically induces a norm continuous shift operator $\textbf{s}$ between two Banach lattices of regular bounded martingales. Moreover, the operators $\hat{T}_{\xi}$ and $\textbf{s}$ commute. Utilizing this fact with inductive limits, we construct a categorical limit space $\mathcal{M}_{T,\xi}$ which is called the associated space of the pair $(T,\xi)$. We present new connections between theories of Boolean algebras, abstract martingales and Banach lattices.

\end{abstract} \maketitle

\section{Introduction}

The notion of martingales in vector lattices goes back to R. DeMarr and J.J. Uhl, see~\cite{DEM} and~\cite{U}. Their algebraic and order theoretical properties are studied by G. Stoica in ~\cite{SG1990,SG1991}. In 2013, G. Stoica used order convergent martingales to investigate market completeness, see~\cite{SG2003}. Theory of martingales in vector lattices has been developed by many researchers since then, see~\cite{SG2006, SG2015, SG2019, EE2012, G, GTX, GT, GL, KLW,T,TX}. In addition, substantial progress has been made in the subject of unbounded convergences on vector lattices, see~\cite{DOT, DE, G, GTX, NG2019} and the references therein. In~\cite{GX}, N. Gao and F. Xanthos used unbounded order Cauchy sequences to generalize Doob's submartingale convergence theorems.

Present work and the results related to $\bb$-Volterra operators on Banach lattices are motivated by both the theory of regular bounded martingales, see~\cite{GX,GT,T,TX,U,SG2019}; and the theory of abstract Volterra operators, see~\cite{FS,Gu,Ku,V,Z,ZA}. 

In the present paper, we define two classes of positive operators on a Banach lattice $E$ with an order continuous norm; namely the class of $\bb$-Volterra operators for a Boolean algebra $\bb$ of order projections, and the class of regular Volterra operators. These operators have the special property that they induce an action on a Banach lattice formed by regular bounded martingales.

Every $\bb$-Volterra operator is a regular Volterra operator. Positive band preserving operators on $E$ provide prime examples of $\bb$-Volterra operators. Conditional expectation operators, see~\cite{KLW}, provide examples for regular Volterra operators. Every regular Volterra operator is a Volterra operator in the sense of~\cite{ZA}. 

Both system and operator theory have developed dramatically to the point where certain classes of operators, see~\cite{AA,FS,FG}, can model systems up to some real world constraints. The class of Volterra operators can be used to produce examples of casual systems in system theory. Much more information related to generalizations of Volterra operators can be found in~\cite{FS,Gu,Ku,V,Z,ZA}. In the settings of system theory, a Volterra operator describes a system whose current state does not depend on future. In the case of abstract martingales in vector lattices, a variant of this property has been already implemented, see~\cite{SG2015, SG2019}. 

Structure of the paper is as follows. In Section~\ref{FirstSec}, we prove several properties of $\bb$-Volterra operators. One of the purposes of this section is to show how lattice operations of vector lattices and Boolean operations of the Boolean algebra $\bb$ can be combined with the properties of a $\bb$-Volterra operator. In Section~\ref{Paths} we define the notion of forward filtrations $\xi$ in $\bb$. The notion of forward filtration can be thought of as a discretization of a resolution of identity in Boolean algebras. In Section~\ref{AbtractMartingaleSpaces}, we define the spaces $\mathcal{M}_0(\xi), \mathcal{M}_{b}(\xi)$ and $ \mathcal{M}_{r}(\xi)$ associated to $\xi$. In Section~\ref{RegularVolterra}, we prove several properties of regular Volterra operators. We show that every regular Volterra operator induces a Volterra operator in the sense of~\cite{ZA}, but not conversely. Discrete case is required for the construction of abstract martingale spaces. We further show in Section~\ref{RegularVolterra} that a conditional expectation operator $T\colon E\to E$ is regular Volterra. For properties of conditional expectation operators see~\cite{KLW}. In Section~\ref{ForwardPaths}, we show that a $\bb$-Volterra operator $T\colon E\to E$ induces a continuous operator $\hat{T}_{\xi}\colon \mathcal{M}_{r}(\xi)\to \mathcal{M}_{r}(\xi)$ on the Banach lattice $\mathcal{M}_{r}(\xi)$ of regular bounded martingales for every forward filtration $\xi$ in $\bb$. By Theorem~\ref{Fundcoro}, the operator $\hat{T}_{\xi}$ commutes with the forward shift operator. Similarly, a regular Volterra operator induces an action on the space of regular bounded martingales for some, note the difference, forward filtration. We then continue with various applications of Theorem~\ref{Fundcoro}. We derive some results about the categorical limit space by \[\mathcal{M}_{T,\xi}\coloneqq \varinjlim (\mathcal{M}_{r}(\xi), \hat{T}_{\xi} )\] which is called the associated space of the pair $(T,\xi)$. In Section~\ref{Antichains}, we show an application related to antichains in the Boolean algebra $\bb$.

\section{Notation and Terminology}\label{Intro}
For unexplained terminology on Boolean algebras, vector lattices and operator theory we refer the reader to~\cite{AB,AA,FG,KK1999,K,PN,S}.

In this paper, $E$ denotes a Banach lattice with an order continuous norm. Classical $L_p(\mu)$-spaces for $1\leq p<\infty$ are prime examples of Banach lattices with an order continuous norm. Detailed discussions on Banach lattices are given in the monographs~\cite{PN,S}, and in~\cite[Chapter 4]{AB},~\cite[par. 1.5.3]{K},~\cite{AA,AA2001}. A Banach lattice is said to ve a \textit{$KB$-space} if every increasing norm bounded sequence in $E$ is norm convergent. If $E$ is a $KB$-space then the norm of $E$ is order continuous, see~\cite{AB,AA}. A \textit{filtration} $(P_n)$ on $E$ is a sequence of positive projections on $E$ such that $P_nP_m=P_mP_n=P_{m\wedge n}$. 

A \textit{Boolean algebra} $\bb$ is a distributive complemented lattice with a distinct zero and one, see~\cite{KK1999,K}. We denote the zero and one of the Boolean algebra $\bb$ by $\textbf{0}$ and $\textbf{1},$ respectively. We denote the complement of $\pi\in \bb$ by $\pi^*$. Evidently, a Boolean subalgebra of $\bb$ always contains the zero and one of $\bb$. A nonempty subset $S\subseteq \bb$ is said to be an \textit{antichain} in $\bb$ if $\pi,\rho\in S$ are distinct elements then $\pi\wedge\rho=0$. 

Denote by $\mathfrak{B}(E)$ the set of all order projections on $E$. The set $\mathfrak{B}(E)$ is a Boolean algebra when it is equipped with the ordering $\pi\leq \rho$ if and only if $\pi \rho =\pi;$ and with lattice operations $\pi\wedge \rho=\pi\rho,$ $\pi\vee \rho=\pi+\rho-\pi\rho$ and $\pi^*=\textbf{1}-\pi$. 

It is a known, see~\cite{AB,AA,KK1999,K,S}, that if $E$ is a Banach lattice with an order continuous norm then there is a one-to-one correspondence between projection bands of $E$ and order projections on $E$. Some extensions of this property can be found in the aforementioned references. In general, the set of ideals of a vector lattice with respect to inclusion is not a complemented lattice, several examples are given in~\cite[Chapter 2]{S}.

Throughout the present paper, $\mathfrak{B}$ denotes a Boolean subalgebra of $\mathfrak{B}(E)$. Therefore, Boolean algebras used in this work are realized by order projections on $E$. An element of $\bb$ uniquely corresponds to an order projection on $E$. We remark that $\pi\rho=\pi\wedge \rho$ and $\pi x=\pi(x)$ for $\pi,\rho\in \bb$ and $x\in E$ because elements of $\bb$ are order projections on $E$. With respect to these conventions, an expression of the form $\pi\wedge \rho x$ is equal to $(\pi\wedge \rho)x$.

Following examples demonstrate how Boolean algebras arise naturally in the study of Banach lattices with order continuous norms. These examples are used in many places in the present article. On the other hand, there are many relationships between Boolean algebras and vector lattices, see the monograph \cite{K} for details.

\begin{exam}\label{SimpleEx} Let $\pi,\rho\in \bb(E)-\{\textbf{0},\textbf{1} \}$ be two order projections. By Theorem 1.45 of~\cite{AB}; $\pi\rho,\pi+\rho$ and $\pi^*$ are order projections on $E$. It follows that the order projections $\pi\rho, \pi\rho^*, \pi, \pi^*\rho, \rho,(\pi\rho^*)\vee(\pi^*\rho), \pi\vee\rho, (\pi\vee\rho)^*, (\pi\rho)\vee(\pi^*\rho^*), \rho^*, \pi\vee\rho^*, \pi^*, \pi^*\vee\rho$ and $(\pi\rho)^*$ belong to $\bb(E)$. It is a corollary of Theorem 1.45 of~\cite{AB} that these elements are all order projections on $E$. Any Boolean subalgebra $\bb$ of $\bb(E)$ containing $\pi$ and $\rho$ also contains all of these order projections.
\end{exam}

We use lowercase letters to denote elements of $E,$ lowercase Greek letters to denote order projections on $E$. This convention is used in many works including~\cite{K}. We write Boolean values $\textbf{0}$ and $\textbf{1}$ in bold with the exception that in some cases $\textbf{1}$ and the identity operator $I_E$ on $E$ are interchangeable. The notation $\mathfrak{B}(E)$ is used by~\cite{K}. We recall that elements of $\mathfrak{B}(E)$ are order projections throughout. 
 
\section{Spaces of $\mathfrak{B}$-Volterra and Regular Volterra Operators }

\subsection{Some Properties of $\mathfrak{B}$-Volterra Operators.}\label{FirstSec}
\begin{definition}\label{Definition00000}
	Let $\bb$ be a Boolean subalgebra of $\bb(E)$. A positive operator $T\colon E\to E$ is said to be \textit{$\mathfrak{B}$-Volterra} if \[ \pi x=\pi y \Rightarrow \pi Tx=\pi Ty \] for every $\pi\in\mathfrak{B}$ and $x,y\in E$. 
\end{definition}

The positivity assumption given in Definiton~\ref{Definition00000} is needed in the sequel.
\begin{prop}
 If a positive operator $T\colon E\to E$ is $\bb$-Volterra then for every $\pi\in \bb$ the vector $\pi Tx$ is uniquely determined by $\pi x$ for every $x\in E$ with $\norm{x}_E=1$.
\end{prop}
\begin{proof}
	Suppose that $\pi x=\pi y$ where $x,y\in E$ satisfy $\norm{x}_E=\norm{y}_E=1$. It follows from Definition~\ref{Definition00000} that $\pi Tx=\pi Ty$.
\end{proof}

A $\bb$-Volterra operator can be compact, and hence a Dunford-Pettis operator. An example is given below. 
\begin{exam}\label{Ex1} Consider the Boolean subalgebra $\bb=\{\textbf{0},\textbf{1} \}$ of $\bb(E)$. In this case, $T\colon E\to E$ is $\bb$-Volterra if and only if $T$ is a positive operator. As another example, let us consider the Boolean subalgebra $\bb=\{\textbf{0},\pi_0,\pi_0^*,\textbf{1} \}$ for some $\pi_0\in \bb(E)$. The operators $\pi_0,\pi_0^*\colon E\to E$ are $\bb$-Volterra. Any nonnegative linear combination $t_1\pi_0+t_2\pi_0^*$ with $t_1,t_2 \geq 0$ is again $\bb$-Volterra. In particular, if $\pi_0$ is a compact order projection then $\pi_0$ is a compact $\bb$-Volterra operator. 
\end{exam}

Following result is motivated from the properties of casual operators on Hilbert resolution spaces, see~\cite{FS} and~\cite[Lemma 10.9]{AA}. Although this result is well-known in the settings of Hilbert spaces, we use it as a tool to deduce new results in the present paper.
\begin{prop}\label{P001}
	Suppose that $T\colon E\to E$ is a positive operator on a Banach lattice $E$ with an order continuous norm. Let $\mathfrak{B}$ be a Boolean subalgebra of $\mathfrak{B}(E)$. The following statements are equivalent:
	\begin{itemize}
		\item[\em i.] {The operator $T$ is $\mathfrak{B}$-Volterra; }
		\item[\em ii.]{$\pi Tx=\pi T \pi x$ for all $x\in E$ and $\pi\in \mathfrak{B};$}
		\item[\em iii.]{$T \pi^* x=\pi^* T \pi^* x$ for all $x\in E$ and $\pi\in \mathfrak{B};$}
		\item[\em iv.]{For every $\pi\in\bb$ the band $\pi(E)$ of $E$ is left fixed by $T$.}
	\end{itemize}
\end{prop}
\begin{proof}
Suppose that $(i)$ holds. Because $\pi^* x\in \ker\pi$ and $T$ is $\mathfrak{B}$-Volterra, $\pi T \pi^* x=\pi T 0=0$. Hence, $\pi T \pi^* x=\pi T (I_E-\pi)x=0$ implies that $\pi Tx=\pi T \pi x$ for all $x\in E$. This shows $(ii)$. Because $\pi Tx=\pi T \pi x$ is satisfied for all $\pi\in \mathfrak{B},$ it is also satisfied for $\textbf{1}-\pi^*\in\mathfrak{B}$. Hence, after cancellations , $(I_E-\pi^*)Tx= (I_E-\pi^*)T(I_E-\pi^*)x$ implies that $\pi^* x=\pi^* T \pi^* x$. The condition given in $(iii)$ implies that the band $\pi(E)$ of $E$ is $T$-invariant. Suppose that $(iv)$ holds. Let $x,y\in E$ be such that $x-y\in\ker\pi$. Hence, $x-y\in (I_E-\pi)(E)$. Because the band $(I_E-\pi)(E)$ is left fixed $T,$ we have $T(x-y)\in (I_E-\pi)(E)$. Hence, $\pi T(x-y)=0$ implies that $\pi Tx=\pi Ty$. Thus, $T$ is a $\mathfrak{B}$-Volterra operator on $E$.
\end{proof}

In view of Proposition~\ref{P001}, we continue with further examples and non-examples of $\bb$-Volterra operators.
\begin{exam}
	Recall from~\cite[Section 4.2]{PN} that a positive $T\colon E\to E$ is said to be band-irreducible if there exists no norm closed $T$-invariant band of $E$ distinct to $\{0\}$ and $E$. We further recall from~\cite[Prop. 5.2]{S} that each band of $E$ is norm closed. Hence, if $\bb$ is a Boolean subalgebra of $\bb(E)$ satisfying $\bb\neq \{ \textbf{0},\textbf{1} \}$ then a band-irreducible $T$ is not $\bb$-Volterra. Because $E$ has an order continuous norm, a similar statement holds for the case of irreducible operators, see~\cite[Section 4.2]{PN}. In details, an irreducible operator can not $\bb$-Volterra.
\end{exam}

\begin{coro}\label{CoroBandPreserving}
	Suppose that $T\colon E\to E$ is a positive operator, and that $\bb$ is a Boolean subalgebra of $\bb(E)$. If $T$ is band preserving then $T$ is $\bb$-Volterra. 
\end{coro}
\begin{proof}
	Result follows from Proposition~\ref{P001}, and from the definition of band preserving operators, see~\cite{AB}. In details, if $T$ is a positive band preserving operator then Proposition~\ref{P001} implies that $T$ is $\bb$-Volterra for every Boolean subalgebra $\bb$ of $\bb(E)$.
\end{proof}

\begin{rem}
	The notions dilatators, orthomorphisms and band preserving operators originate from the works of H. Nakano, A. Bigard, K. Keimel, P.F. Conrad, J.E. Diem and M. Meijer. A detailed discussion in this direction can be found in~\cite[page 659]{Z1983}.
\end{rem}

\begin{prop}\label{P002}
	Suppose that $T\colon E\to E$ is a $\mathfrak{B}$-Volterra operator on an order continuous Banach lattice $E$ for some Boolean subalgebra $\mathfrak{B}$ of $\mathfrak{B}(E)$. The positive operator $T^n$ is a $\mathfrak{B}$-Volterra operator for all $n\geq 0$. Composition of two $\mathfrak{B}$-Volterra operators is again $\mathfrak{B}$-Volterra. Further, the set of $\bb$-Volterra operators generates a band in the ordered vector space $\mathcal{L}_r(E)$ of all regular operators from $E$ to $E$.
\end{prop}
\begin{proof} Let $\pi\in\bb$. By Proposition~\ref{P001}, if $T$ is $\mathfrak{B}$-Volterra then the band $\pi(E)$ of $E$ is left fixed by $T,$ and hence $\pi(E)$ is left fixed by $T^n$ for all $n\geq 1$. It follows that $T^n$ is $\mathfrak{B}$-Volterra for all $n\geq 1$. In the case $n=0,$ the operator $T^0=I_E$ is again $\bb$-Volterra. This follows from Definition~\ref{Definition00000}. 

Let $T_1,T_2$ be two $\mathfrak{B}$-Volterra operators. The operators $T_1T_2$ and $(T_1+T_2)^2$ are positive operators on $E$. Evidently $T_1T_2$ and $(T_1+T_2)^2$ are again $\mathfrak{B}$-Volterra.

Because $\bb$-Volterra operators are positive operators, and positive operators are regular by~\cite[page 12]{AB}; there exists a unique smallest band in the ordered vector space $\mathcal{L}_r(E)$ containing all $\bb$-Volterra operators on $E$.
\end{proof}

It is stated in Proposition~\ref{P002} that collection of $\bb$-Volterra operators generates a band, denoted by $L_{\bb}(E)$, in the space $\mathcal{L}_r(E)$ of all regular operators on $E$. We remark that not all operators belonging to $L_{\mathfrak{B}}(E)$ are $\bb$-Volterra. 

We further note that the assignment $\bb\mapsto L_{\bb}(E)$ is contravariant in the sense that if $\bb_1$ is a Boolean subalgebra of $\bb_2$ then $L_{\bb_2}(E)$ is a subspace of $L_{\bb_1}(E)$. By using the regular norm $||.||_r$ on $\mathcal{L}_r(E),$ see~\cite[Chapter 4.4]{AB}, we can write \[\norm{T}_r=\inf \{\norm{S}\colon T\leq S \mbox{ where } S\colon E\to E \mbox{ is positive}\} \] for a $\bb$-Volterra operator $T$ where $\norm{S}$ denotes the operator norm of $S$. Because $T$ is positive, $\norm{T}=\norm{T}_r$. Further, as every band of a vector lattice is a vector lattice in itself, $(L_{\bb}(E),\norm{\cdot}_r)$ forms a normed vector lattice.

\begin{lem}\label{Lemma1}
	Suppose that $T\colon E\to E$ is a positive operator. Let $\mathfrak{B}$ be a Boolean subalgebra of $\mathfrak{B}(E)$. 
	\begin{itemize}
		\item[\em i.] { If $T$ is $\bb$-Volterra then the operators $\pi T\colon E\to E$ and $\pi+T-\pi T\colon E\to E$ are $\bb$-Volterra for every $\pi\in \bb$. Further, $(\pi T)^k=\pi T^k$ for all $k\geq 1$. }
		\item[\em ii.]{ Let $A$ be an ideal in $E,$ and suppose that $\bb$ is the Boolean algebra generated by order projections $\pi_{x}(z)\coloneqq \sup\{z\wedge n|x|\colon n\in\nn\}$ for $x\in A,$ and $z\in E$. If $T$ is $\bb$-Volterra then $T$ is $\bb'$-Volterra for every Boolean subalgebra $\bb'$ generated by $\{\pi_x\colon x\in A' \}$ where $A'\subseteq A$ is nonempty. }
		\item[\em iii.] { If $T$ is $\bb$-Volterra then $\pi T(t\pi x+(1-t)x)=\pi Tx$ for all $t\in [0,1],$ $x\in E$ and $\pi\in\bb$. }
		\item[\em iv.] { For every $\pi \in \bb$ denote by $I_{\pi}$ the linear subspace of $L_{\bb}(E),$ see Proposition~\ref{P002}, generated by those $S\in L_{\bb}(E)$ satisfying $\pi S=S\pi$. If $T$ is $\bb$-Volterra then any linear combination of operators of the form $t_1\pi T+t_2\pi^*T+t_3\rho T$ belongs to $I_{\pi}$ where $t_1,t_2,t_3\geq 0$ and $\rho\in \bb$ satisfies $\pi\wedge \rho=\textbf{0}$. }
	\end{itemize}
\end{lem}
\begin{proof}
$(i)$. We first observe that because $\pi\in \bb$ is an order projection on $E$, $\pi T$ is a positive operator. Let $\rho\in \bb$ be arbitrary. For $x,y\in E,$ the equality $\rho x=\rho y$ implies that $(\rho \wedge \pi)\rho x=(\rho \wedge \pi)\rho y$. Because $\rho \wedge \pi\leq \rho$ in $\bb$, we have $(\rho \wedge \pi)\rho=\rho \wedge \pi$. Hence, we obtain $(\rho \wedge \pi) x=(\rho\wedge \pi) y$. The operator $T$ is $\bb$-Volterra. We obtain \[(\rho \wedge \pi) Tx=\rho \pi Tx=\rho \pi Ty=(\rho \wedge \pi) Ty.\] The last equality shows that $\pi T$ is $\bb$-Volterra. 

In order to show that $\pi+T-\pi T\colon E\to E$ is $\bb$-Volterra, we observe that $\pi+T-\pi T=\pi+(I_E-\pi)T=\pi +\pi^*T$. Hence, by Proposition~\ref{P002} and the above paragraph, $\pi+T-\pi T$ is $\bb$-Volterra. 

$(ii)$. We first show the simple fact that if a nonempty subset $B\subseteq \bb$ generates the Boolean subalgebra $\bb'$ of $\bb$ then $T$ is $\bb'$-Volterra. Because the operator $T$ is $\bb$-Volterra, $T$ is $\bb'$-Volterra for every Boolean subalgebra $\bb'$ of $\bb$. In particular, this is the case if $\bb'$ is generated by the nonempty subset $B\subseteq \bb$.

Since $E$ is order complete, every vector belonging to $E$ is a projection vector in $E$. Let $x\in A$. By~\cite[Theorem 1.47]{AB} and~\cite[Corollary 2.11.2]{S}, the formula $\pi_{x}(z)\coloneqq \sup\{z\wedge n|x|\colon n\in\nn\}$ for $z\in E$ defines an order projection, and hence $\pi_x$ belongs to $\bb(E)$. The Boolean algebra $\bb$ is the smallest Boolean subalgebra of $\bb(E)$ containing the collection $\{ \pi_x\colon x\in A\}\subseteq \bb(E)$. From the first part of the present proof, conclusion follows. 

$(iii)$. Let $\pi\in\bb$ and $t\in [0,1]$. It follows from Proposition~\ref{P001} that \[\pi T(t\pi x+(1-t)x)=t\pi T\pi x+(1-t)\pi Tx=\pi Tx \] for all $x\in E$.

$(iv)$. Let $\pi\in \bb$ be fixed. Consider the operator $S=t_1\pi T+t_2\pi^*T+t_3\rho T$ for $t_1,t_2,t_3\geq 0$ and $\pi\wedge \rho=\textbf{0},$ $\rho\in \bb$. Because we have \[ \pi S=t_1\pi T + t_2\pi\pi^* T +t_3\pi\rho T=t_1 \pi T \] and \[S\pi=t_1\pi T\pi+t_2\pi^*T\pi+t_3\rho T \pi=t_1\pi T+t_2\pi^*\pi T\pi+t_3\rho \pi T \pi=t_1 \pi T \] by Proposition~\ref{P001}, we have $\pi S= S\pi$. The operators $\pi T,$ $\pi^*T$ and $\rho T$ are $\bb$-Volterra for every $\pi,\rho\in\bb$. Hence, $S=t_1\pi T+t_2\pi^*T+t_3\rho T$ is $\bb$-Volterra by Proposition~\ref{P002}. It follows that $S\in I_{\pi}$. A linear combination of such operators commutes with $\pi,$ and hence it belongs to $I_{\pi}$. Indeed, if $S'=\sum_{j=0}^ns_jS_j$ with $S_j=t_{1,j}\pi T_j+t_{2,j}\pi^*T_j+t_{3,j}\rho_j T_j$ for some $\bb$-Volterra operators $T_1,T_2,\ldots,T_n;$ $t_{1,j},t_{2,j},t_{3,j}\geq 0$ and $\pi\wedge \rho_j=\textbf{0},$ then $\pi S'=S'\pi$.
\end{proof}

\begin{prop} Suppose that $T\colon E\to E$ is a $\mathfrak{B}$-Volterra operator on an order continuous Banach lattice $E$ for some Boolean subalgebra $\mathfrak{B}$ of $\mathfrak{B}(E)$. Suppose further that $0\leq x\leq y\leq z$ in $E,$ and that $\pi x\wedge \rho y=\pi y\wedge \rho z$ holds for some $\pi,\rho\in \bb$. Then $\pi\rho T^k(t\pi x+(1-t)x)=\pi \rho T^ky$ for every $k\geq 1$ and $t\in [0,1]$.
\end{prop}
\begin{proof} We use the fact that if $\pi$ and $\rho$ are two order projections then $\pi x\wedge \rho y= \pi\rho (x\wedge y)$ for all $x,y\in E^+$ see~\cite[Section 1.3]{AB}. Let $0\leq x\leq y\leq z$ and $\pi,\rho$ be as given in hypothesis. It follows that \[\pi x\wedge \rho y=\pi\rho(x\wedge y)=\pi\rho x=\pi\rho y=\pi\rho(y\wedge z)=\pi y\wedge \rho z. \] Because $T$ is $\bb$-Volterra, we conclude that $\pi\rho T^kx=\pi \rho T^ky$ for every $k\geq 1$. As a result of Lemma~\ref{Lemma1}, we further have $\pi\rho T^k(t\pi x+(1-t)x)=\pi \rho T^ky$ for every $t\in [0,1]$ and $k\geq 1$.
\end{proof}
The operator $\pi-\pi T + \pi^* T$ given in Theorem~\ref{P0010} is abstracted from the order projection $(\pi\rho^*)\vee(\pi^*\rho),$ see Example~\ref{SimpleEx}. 
\begin{thm}\label{P0010} Suppose that $T\colon E\to E$ is a $\mathfrak{B}$-Volterra operator satisfying $ T\leq I_E$ on an order continuous Banach lattice $E$ for some Boolean subalgebra $\mathfrak{B}$ of $\mathfrak{B}(E)$. Then the operator $\pi-\pi T + \pi^* T\colon E\to E$ is $\bb$-Volterra for every $\pi\in \bb$. Further, if $\pi\in \bb$ satisfies $\pi T=\pi$ then the operator $I_E+\pi-T$ is $\bb$-Volterra.
\end{thm}
\begin{proof} We first show that $\pi-\pi T+\pi^*T$ is a positive operator. Because $T$ is positive and $T\leq I_E$ we have $0\leq \pi (I_E-T)$. Hence, $\pi T\leq \pi$. This implies that $-\pi T\leq \pi-2\pi T$ and that $\pi^*T=T-\pi T\leq \pi-2\pi T+T$. Hence, $\pi-2\pi T+T\geq \pi^* T$. The operators $\pi^*$ and $T$ are positive operators. Hence, $\pi-\pi T+\pi^*T=\pi-2\pi T+T\geq 0$.

As observed previously $\pi-\pi T+\pi^*T=\pi-2\pi T+T$. Suppose $\rho x=\rho y$ for some $\rho\in \bb$ and $x,y\in E$. Because $T$ is $\bb$-Volterra, and, $\pi$ and $\rho$ commute; it follows that \[ \rho \pi x-2\pi \rho Tx+\rho Tx=\pi\rho y-2\pi \rho Ty+\rho Ty .\] As $\rho\in\bb$ is arbitrary, the last equality shows that $\pi-\pi T+\pi^*T$ is $\bb$-Volterra.

For the second part, observe that $\pi T=T$ implies \[I_E+\pi-T=\pi^*+2\pi-T=\pi^*+2\pi T-T.\] We first show that $\pi^*+2\pi T-T$ is positive whenever $\pi T=T$ holds. We have $ \pi+T\leq I_E+\pi T$ because by Theorem 1.44 of~\cite{AB} $\pi\leq I_E$. Hence, $\pi T\leq I_E-\pi+2\pi T-T$. This implies that $\pi T\leq \pi^*+2\pi T-T$. Both $\pi$ and $T$ are positive operators. Hence, $\pi^*+2\pi T-T\geq 0$. Finally, suppose $\rho x=\rho y$ for some $\rho\in \bb$ and $x,y\in E$. It follows that \[\rho\pi^*x+2\pi\rho Tx-\rho Tx= \rho \pi^*y+2\pi\rho Ty-\rho T y \] because $T$ is $\bb$-Volterra. Hence, $\pi^*+2\pi T-T$ is $\bb$-Volterra.
\end{proof}

\begin{prop} Let $A$ be an ideal in $E,$ and denote by $\bb$ the Boolean subalgebra of $\bb(E)$ generated by order projections of the form $\pi_{x}(z)\coloneqq \sup\{z\wedge n|x|\colon n\in\nn\}$ for $x\in A$ and $z\in E$. Suppose that $T$ is $\bb$-Volterra, $0\leq x_3\leq x_2\leq x_1$ in $A,$ and that $\{\pi_{x_1},\sigma,\rho \}$ is an antichain in $\bb$. Then \[ (\pi_{x_1}-\pi_{x_2})T\sigma x\perp (\pi_{x_2}-\pi_{x_3})T(\rho^*\sigma x+w) \] for every $x\in E$ and $w\in \pi_{x_1}^*(E)$. In particular, $(\pi_{x_1}-\pi_{x_2})T\sigma x \perp \pi_{x_2}T\rho^*\sigma x$ for every $x\in E$.
\end{prop}
\begin{proof} We first observe that \[ [(\pi_{x_1}-\pi_{x_2})T\sigma x]\wedge [(\pi_{x_2}-\pi_{x_3})T(\rho^*\sigma x+w)] \] 
	is equal to
	\[ [(\pi_{x_1}-\pi_{x_2})\pi_{x_1}T\sigma x]\wedge [(\pi_{x_2}-\pi_{x_3})\pi_{x_1}T(\rho^*\sigma x+w)] \] and that $\pi_{x_1}T(\rho^*\sigma x+w)=\pi_{x_1}T\rho^*\sigma x$ for every $w\in \pi_{x_1}^*(E)$.
	
	Let us show that $\pi_{x_1}T\sigma x=\pi_{x_1}T(\rho^*\sigma x)$ for every $x\in E$. Because $\{\pi_{x_1},\sigma,\rho \}$ is an antichain in $\bb,$ we have $\pi_{x_1}\wedge \sigma=\textbf{0}$ and $\rho\wedge \sigma=\textbf{0}$. Hence, \[\pi_{x_1} T\rho\sigma x=\pi_{x_1} T(I_E-\rho^*)\sigma x=0\] for all $x\in E$. Thus, $\pi_{x_1}T\sigma x=\pi_{x_1}T(\rho^*\sigma x)$ for every $x\in E$.
	
	If we put $z=\pi_{x_1}T\sigma x$ then we have $(\pi_{x_1}-\pi_{x_2})z\wedge (\pi_{x_2}-\pi_{x_3}) z=0,$ see~\cite[Theorem 1.48]{AB}. Particular case follows by taking $x_3=0$.
\end{proof}

\subsection{Forward and Backward Filtrations in $\bb$}\label{Paths}
The notion of forward filtration is needed for the construction of the abstract martingale spaces given in the sequel. It can be thought of as a discretization of a resolution of identity in Boolean algebras.
\begin{definition}
A function of the form \[\xi\colon \hat{\mathbb{N}}=\{0,1,\ldots,\infty\}\to \bb\] is said to be a \textit{forward filtration} in $\bb$ if 
\begin{enumerate}
	\item $\xi_n\coloneqq\xi(n)\leq \xi_{n+1}$ for every $n\geq 0$ where $\leq$ is the lattice order of $\bb$,
	\item $\xi_0=\bold{0}$ and $\xi_{\infty}=\bold{1}$. 
\end{enumerate}
\end{definition}

We denote by $\mathcal{P}^+(\bb)$ the set of all forward filtrations in $\bb$. In some cases, it is convenient to use the tabular form \[\xi:\xi_0=\textbf{0}\leq \xi_1\leq \xi_2 \leq \cdots \leq \xi_{\infty}=\textbf{1}\] to denote the forward filtration $\xi$.

In order to construct a dual to a forward filtration, we use the natural idempotent map $^*\colon \bb\to \bb$ of Boolean algebras. Because filtrations are formed by increasing sequence of projections, the dual notion with respect to the idempotent map $^*\colon \bb\to \bb$ may not be a filtration itself, also see~\cite[Sec.11]{T}.

\begin{definition}
	A function $\xi\colon \hat{\nn}\to \bb$ is said to be a \textit{backward filtration} in $\bb$ if 
	\begin{enumerate}
		\item $\xi_{n+1}\coloneqq \xi(n+1)\leq \xi_{n}$ for every $n\geq 0$ where $\leq$ is the lattice order of $\bb$,
		\item $\xi_0=\bold{1}$ and $\xi_{\infty}=\bold{0}$.
	\end{enumerate}
\end{definition}

We denote by $\mathcal{P}^-(\bb)$ the set of all backward filtrations in $\bb$. The backward filtration $\xi$ can be written in the tabular form as \[\xi\colon \xi_0=\textbf{1}\geq \xi_1\geq \xi_2\geq \cdots\geq \xi_{\infty}=\textbf{0}.\] It follows that if $\xi\in \mathcal{P}^-(\bb)$ then $(\xi^*)_n\coloneqq (\xi_n)^*$ for $n\in \hat{\nn}$ defines a forward filtration, that is $\xi^*\in \mathcal{P}^+(\bb)$.

\begin{rem}\label{Remark1}
	In~\cite{GX}, the definition of abstract bistochastic filtrations can be found. An abstract bistochastic filtration assumes the existence of a weak unit in $E$ which is left fixed by all projections belonging to filtration. Evidently, the notions of abstact bistochastic and forward filtrations are nonequivalent in the settings of order continuous Banach lattices and order projections. We further note that a more general approach to filtrations can be found in~\cite{SG1990, SG1991}.
\end{rem}

In view of Remark~\ref{Remark1}, concrete and natural examples of forward filtrations can be obtained by discretizing resolutions of the identity in Boolean algebras. For the definition of resolution of the identity in Boolean algebras, see~\cite[1.4.3]{K}. 
\begin{thm}\label{TrivialProp}
	Let $e\colon \rr\to \bb$ be a resolution of the identity such that $e(0)=\textbf{0}$. For every strictly increasing and divergent sequence $t_n$ in $\mathbb{R}$ the function $\xi\colon\hat{\nn}\to \bb$ defined by $\xi(n)\coloneqq e(t_n)$ for $n\in \nn$ and $\xi(\infty)\coloneqq\textbf{1}$ is a forward filtration in $\bb$. Further, if $e_0\colon \rr\to \bb$ is a resolution of the identity satisfying $e_0(s_0)=\textbf{0}$ for some $s_0\in \rr$ then $e(t)\coloneqq e_0(t+s_0)$ is a resolution of the identity satisfying $e(0)=\textbf{0}$.
\end{thm}
\begin{proof} If $s\leq t$ then $e(s)\leq e(t)$. Because $e(0)=\textbf{0},$ we have $e(t)=\textbf{0}$ for all $t\leq 0$. Since $\bigvee_{t\in \rr} e(t)=\textbf{1},$ there is a sequence in $[0,\infty]\subseteq \rr$ such that $e(t_n)\leq e(t_{n+1})$ for all $n$. Hence, $\xi$ is a forward filtration in $\bb$. For the second part, we observe that if $e_0\colon \rr\to \bb$ is a resolution of the identity then $e(t)\coloneqq e_0(t+s_0)$ is again a resolution of the identity. Evidently $e_0(s_0)=\textbf{0}$ implies that $e(0)=\textbf{0}$. 
\end{proof}

Following result shows how to transform a pairwise order disjoint sum in a Banach lattice into a forward filtration. This transformation is well-defined up to an ordering of summands in pairwise order disjoint sums.
\begin{coro}
	Let $E$ be a Banach lattice with an order continuous norm. 	For every finite partition of unity consisting of nonzero order projections $\pi_1,\pi_2,\ldots,\pi_n$ in $\bb$ the assignment $\xi_n\coloneqq \bigvee_{k=1}^n\pi_k$ with $\xi_k=\textbf{1}$ for $k\geq n+1$ together with $\xi_0=\textbf{0}$ and $\xi_{\infty}=\textbf{1}$ define a forward filtration in $\bb$. For every pairwise order disjoint sum $x_1+x_2+\dots+x_n$ in $E$ there is a corresponding forward filtration in $\bb(E)$.
\end{coro}
\begin{proof}
If $\pi_1,\pi_2,\ldots,\pi_n$ is a finite partition of unity then the given formula guarantees that $\xi$ is a forward filtration because $\xi_n\leq \xi_{n+1}$ for all $n\in \hat{\nn}$. 
For every pairwise order disjoint sum $x_1+x_2+\dots+x_n$ there exists a finite partition unity $\pi_1,\pi_2,\ldots,\pi_n$ in $\bb(E)$ such that $\pi_kx_l=0$ if $k\neq l$ where $k=1,2,\ldots,n;$ see~\cite[2.1.9]{K} and~\cite[Lemma 2.2]{AB}. 
\end{proof}

We note that there is a certain symmetry of the set $\mathcal{P}^+(\bb)$, i.e., there is a discrete semigroup acting on $\mathcal{P}^+(\bb)$. In details, for every $\pi\in \bb$ let us define \[\mbox{Sc}(\pi)\coloneqq \{\pi'\in\bb\colon \pi\leq\pi'\}.\] In view of~\cite[Chapter 2]{S}, we say that $\mbox{Sc}(\pi)$ is a section. We further say that a transformation $S\colon \bb\to\bb$ is \textit{sectionally open} if $\mbox{Sc}$ and $S$ commute. Semigroup of sectionally open transformations fixing $\textbf{0}$ reappears in Theorem~\ref{StarOpenInduced} of Section~\ref{ForwardPaths}. In general, it is easy to find a sectionally open transformation $S\colon \bb\to\bb$ not fixing $\textbf{0}$.
\begin{prop}\label{StarOpen}
	Let $S\colon \bb\to\bb$ be a sectionally open transformation fixing $\textbf{0}\in \bb$. If $\xi\in \mathcal{P}^+(\bb)$ then the map $n\mapsto S(\xi(n))=S(\xi_n)$ for $n\in\hat{\nn}$ belongs to $\mathcal{P}^+(\bb)$. In particular, the discrete semigroup of sectionally open transformations fixing $\textbf{0}$ acts on $\mathcal{P}^+(\bb)$ via coordinatewise evaluation.
\end{prop}
\begin{proof}
	Because $\xi\in\mathcal{P}^+(\bb) $ implies $\xi_0=\textbf{0},$ it suffices to show that $S(\xi_n)\leq S(\xi_{n+1})$ for all $n\geq 1$. The inequality $\xi_n\leq \xi_{n+1}$ implies that $\xi_{n+1}\in \mbox{Sc}(\xi_n)$. Therefore, $S(\xi_{n+1})\in S(\mbox{Sc}(\xi_n))=\mbox{Sc}(S(\xi_n))$. This means that $S(\xi_{n+1})\geq S(\xi_{n})$. Evidently $S(\textbf{1})=\textbf{1}$ whenever $S$ is a sectionally open transformation. This implies that $S(\xi_{\infty})=\textbf{1}$.
\end{proof}
\subsection{Regular and Bounded Martingales.}\label{AbtractMartingaleSpaces}
For a forward filtration $\xi \in \mathcal{P}^+(\bb),$ we define \[ \mathcal{M}_0(\xi)\coloneqq \{(x_n)_{n=1}^\infty\colon \xi_nx_m=x_n \:\:\mbox{if $m\geq n$} \}\] where $(x_n)_{n=1}^\infty$ is a sequence in $E$. We remark the obvious but important fact that $n$ appearing in the definition of $\mathcal{M}_0(\xi)$ is a finite index varying from 1 to $\infty$. For a backward filtration $\xi \in \mathcal{P}^-(\bb),$ we define \[ \mathcal{M}_1(\xi)\coloneqq \{(x_n)_{n=1}^\infty\colon (I-\xi_n)x_m=x_n \:\:\mbox{if $m\geq n$} \}.\]

The following result reflects one of the main ideas of the present article. In the terminology of ~\cite{SG1990, SG1991}, we have $(x_n)_{n=1}^\infty\in \mathcal{M}_0(\xi)$ if and only if $(x_n)_{n=1}^\infty$ is $(\xi_n(E))_{n=1}$-adapted.

\begin{prop}
	If $\xi\in \mathcal{P}^+(\mathfrak{B})$ then the sequence $(\xi_n)_{n=1}^{\infty},$ note that $\xi_0=\textbf{0}$ and $\xi_{\infty}=\textbf{1}$ are excluded, is a filtration on $E$ in the sense of~\cite{GX,T,TX,U}. The space $\mathcal{M}_0(\xi)$ is identified with the space of all abstract martingales with respect to this filtration. The identification map is tautological.
\end{prop}

Following~\cite{T,TX}, given $\xi\in \mathcal{P}^+(\mathfrak{B})$ we put 
\[ \mathcal{M}_{b}(\xi)\coloneqq \{(x_n)_{n=1}^\infty\in \mathcal{M}_{0}(\xi)\colon \sup_{n\geq 1}||x_n||<\infty \}\subseteq (\bigoplus_{n=1}^{\infty}E)_{\infty} .\] The space $\mathcal{M}_{b}(\xi)$ is called \textit{the space of all bounded martingales on} $E$ corresponding to the forward filtration $\xi\in \mathcal{P}^+(\mathfrak{B})$. We further put
\[ \mathcal{M}_{r}(\xi)\coloneqq \{x^1-x^2\colon x^{i}\in (\mathcal{M}_{b}(\xi))^+ \mbox{ for } i=1,2 \}\] for each $\xi\in \mathcal{P}^+(\mathfrak{B})$. The space $\mathcal{M}_{r}(\xi)$ is called \textit{the space of all regular bounded martingales} corresponding to the forward filtration $\xi\in \mathcal{P}^+(\mathfrak{B})$. 
The regular martingale norm on $ \mathcal{M}_{r}(\xi)$ is given by
\[ \norm{x}_r=\inf_{\substack{y\in(\mathcal{M}_{b}(\xi))^+\\ \pm x\leq y}} \sup_{n\geq 1}\norm{y_n} \] for $x\in \mathcal{M}_{r}(\xi)$. 

\begin{exam}
	Let us give two examples of abstract martingales; they also illustrate how the associated abstract martingale spaces varies as the forward filtration varies. First, consider the forward filtration $\xi_n=\textbf{0}$ for $0\leq n<\infty$ and $\xi_{\infty}=\textbf{1}.$ It follows that both $\mathcal{M}(\xi)$ and $\mathcal{M}_b(\xi)$ are trivial spaces. Second, consider the forward filtration $\xi_0=\textbf{0}$ and $\xi_{n}=\textbf{1}$ for $1\leq n\leq \infty$. In this case, $\mathcal{M}(\xi)=\prod_{n=1}^{\infty}E$ and $\mathcal{M}_b(\xi)=(\bigoplus_{n=1}^{\infty}E)_{\infty}$. 
\end{exam}

\subsection{Regular Volterra Operators}\label{RegularVolterra}
In Section~\ref{FirstSec}, we obtained some results related to $\bb$-Volterra operators. Our present focus is the space of regular Volterra operators. 

A regular Volterra operator can be considered as a Volterra operator in the sense of~\cite{ZA} and of~\cite{Gu,Ku,V,Z}. Furthermore, we show that a conditional expectation operator, see~\cite{KLW}, is a regular Volterra operator.

\begin{definition}
We say that a positive operator $T\colon E\to E$ is a \textit{regular Volterra operator} with respect to a forward filtration $\xi\in \mathcal{P}^+(\bb(E))$ if \[ \xi_n x=\xi_n y\Rightarrow \xi_n Tx=\xi_n T y\]
for all $n\in\hat{\nn}$ and $x,y\in E$.
\end{definition}

\begin{exam}
	This example is motivated from~\cite[Example 1.7]{V}. More information about Schauder bases and complemented subspaces of order continuous Banach lattices can be found in~\cite[Section 5.1]{PN}. Let $E$ be an order continuous Banach lattice with a Schauder base $(e_i)_{i=1}^{\infty}$ such that the canonical projections $P_n\colon E\to E$ defined by $P_n(\sum_{i=1}^{\infty}a_ie_i)=\sum_{i=1}^{n}a_ie_i$ are order projections. Suppose that \[ T(\sum_{i=1}^{\infty} a_ie_i)=\sum_{i=1}^{\infty} (\sum_{k=1}^{\infty}g_{ik}(a_k))e_i\] with $g_{ik}(a_k)=0$ for $k>i$. Then $T$ is a positive Uryson operator. It follows that $T$ is regular Volterra with respect to forward filtration given by the tabular form $\xi: P_0=\textbf{0}\leq P_1\leq P_2\leq \cdots \leq P_{\infty}=\textbf{1} $.
\end{exam}

In the next result, we obtain new positive operators by using a regular Volterra operator. These new positive operators commute with certain projections.
\begin{prop}
	Let $T\colon E\to E$ be a regular Volterra operator with respect to some forward filtration $\xi$. Denote by $\xi_{n_1}\leq \xi_{n_2}\leq \cdots \leq \xi_{n_k}$ the finite chain of order projections whose indices form a finite chain $n_1\leq n_2\leq \cdots \leq n_k$ in $\hat{\nn}$. Consider the subspace \[ B_{n_1\leq n_2\leq \cdots \leq n_k} \coloneqq \{(x_1,x_2,\ldots,x_k)\colon \xi_{n_i}(x_i)=\xi_{n_j}(x_j) \mbox{ for all } i,j\in \{ 1,2,\ldots,k\} \}\] of the product $\prod_{i=1}^k E$.
	\begin{itemize}
		\item[\em i.] { The operator $T$ induces an additive, positively homogeneous operator \[T_{n_1\leq n_2\leq \cdots \leq n_k}\colon B_{n_1\leq n_2\leq \cdots \leq n_k} \to B_{n_1\leq n_2\leq \cdots \leq n_k}.\] }
		\item[\em ii.] { The operator $T$ induces an additive, positively homogeneous operator \[T_{F}'\colon \sum_F B_{n_1\leq n_2\leq \cdots \leq n_k} \to E\] for every finite set $F$ of finite chains $n_1\leq n_2\leq \cdots \leq n_k$. }
		\item[\em iii.] { Given two chains of the form $n_1\leq n_2\leq \cdots \leq n_{k+1}$ and $n_1\leq n_2\leq \cdots \leq n_k$ for some $k\geq 1,$ there exists a positive canonical projection \[\pi\colon B_{n_1\leq n_2\leq \cdots \leq n_{k+1}}\to B_{n_1\leq n_2\leq \cdots \leq n_k}\] such that $\pi\circ T_{n_1\leq n_2\leq \cdots \leq n_{k+1}}=T_{n_1\leq n_2\leq \cdots \leq n_{k}} \circ \pi$. }
	\end{itemize} 
\end{prop}
\begin{proof}
$(i)$. Let us define $T_{n_1\leq n_2\leq \cdots \leq n_k}(x_1,x_2,\ldots,x_k)=(Tx_1,Tx_2,\ldots, Tx_k)$. Because $T$ is a regular Volterra operator with respect to $\xi,$ we have $\xi_{n_i}x=\xi_{n_i}y$ implies $\xi_{n_i} Tx=\xi_{n_i} Ty$ for all $1\leq i\leq k$. This shows that $(Tx_1,Tx_2,\ldots, Tx_k)\in B_{n_1\leq n_2\leq \cdots \leq n_k}$ whenever $(x_1,x_2,\ldots,x_k)\in B_{n_1\leq n_2\leq \cdots \leq n_k}$. If $(x_1,x_2,\ldots,x_k)$ is positive with respect to coordinatewise order then $(Tx_1,Tx_2,\ldots, Tx_k)$ is positive. 

$(ii)$. For every finite chain $n_1\leq n_2\leq \cdots \leq n_{k}$ there exists a projection $\pi_{n_1\leq n_2\leq \cdots \leq n_{k}}\colon B_{n_1\leq n_2\leq \cdots \leq n_k}\to E$ defined by $(x_1,x_2,\ldots,x_k)\to x_1$. Hence, \[\pi_{n_1\leq n_2\leq \cdots \leq n_{k}}\circ T_{n_1\leq n_2\leq \cdots \leq n_k}\colon B_{n_1\leq n_2\leq \cdots \leq n_k}\to E\] induces $T_{F}'\colon \sum_F B_{n_1\leq n_2\leq \cdots \leq n_k} \to E$ by \[T_F'=\sum_F \pi_{n_1\leq n_2\leq \cdots \leq n_{k}}\circ T_{n_1\leq n_2\leq \cdots \leq n_k}\] where the sum runs through all chains $n_1\leq n_2\leq \cdots \leq n_k$ belonging to $F$.

$(iii)$. If $(x_1,x_2,\ldots,x_{k+1})\in B_{n_1\leq n_2\leq \cdots \leq n_{k+1}}$ then $(x_1,x_2,\ldots,x_{k})\in B_{n_1\leq n_2\leq \cdots \leq n_{k}}$. Hence, consider $\pi(x_1,x_2,\ldots,x_{k+1})=(x_1,x_2,\ldots,x_k)$. It follows that \[ \pi \circ T_{n_1\leq \cdots \leq n_{k+1}} (x_1,\ldots,x_{k+1})=(Tx_1,\ldots,Tx_k)=T_{n_1\leq \cdots \leq n_{k}} \circ \pi (x_1,\ldots,x_{k+1})\] for every $(x_1,\ldots,x_{k+1})\in B_{n_1\leq n_2\leq \cdots \leq n_{k+1}}$.
\end{proof}

It follows, details are given below, that a regular Volterra operator can be considered as a Volterra operator in the sense of~\cite{ZA} and of~\cite{Gu,Ku,V,Z}. However, it is possible to have a Volterra operator on $E$ in the sense of~\cite{ZA} which is not a regular Volterra operator, in general. 

\begin{prop}\label{Transfer1} Let $T\colon E\to E$ be a regular Volterra operator on an order continuous Banach lattice $E$. Then there is a system, in the sense of~\cite{ZA}, on $E$ with respect to which $T$ is a Volterra operator in the sense of~\cite{ZA}.
\end{prop}
\begin{proof}
Let $\phi\colon (0,\infty)\to (0,1)$ be an order reversing homeomorphism of the real interval $(0,\infty)$ onto $(0,1)$. We define relations \[xR_t y \mbox{ \: for $x,y\in E$ if and only if \:} x-y\in \xi_{\lfloor \phi^{-1}(t) \rfloor}(E) \] for $t\in (0,1)$ together with $R_0=E\times E$ and $R_1=\Delta_E,$ the diagonal of $E,$ then the one-parameter set of equivalence relations $\mathcal{B}=\{R_t\colon t\in [0,1] \}$ becomes a system in the sense of~\cite{ZA}. Because $T$ is a regular Volterra operator with respect to $\xi,$ the operator $T$ is a Volterra operator on the system $\mathcal{B}$ in the sense of~\cite{ZA}.
\end{proof}

We present two different proofs of the following result. First proof uses vector lattice properties of $E$. Second proof uses contraction of the one-parameter set of equivalence relations, see~\cite[Sec. 2]{ZA}, appearing in the proof of Proposition~\ref{Transfer1}.
\begin{prop}
	Let $T\colon E\to E$ be a regular Volterra operator with respect to $\xi\in P^+(\bb(E))$. Let $k\in\hat{\nn}$ be fixed. The restricted operator $T|_{\xi_k(E)}\colon \xi_k(E)\to \xi_k(E)$ is a regular Volterra operator with respect to some forward filtration in the Boolean algebra $\bb(\xi_k(E))$.
\end{prop}
\begin{proof} As $\xi\in P^+(\bb(E)),$ the function $\xi'\colon \hat{\nn}\to \bb(\xi_k(E))$ defined by $\xi'(m)=\xi_m\xi_k$ belongs to the set $P^+(\xi_k(E))$. Evidently the operator $T|_{\xi_k(E)}\colon \xi_k(E)\to \xi_k(E)$ is positive on the Banach sublattice $\xi_k(E)$ of $E$. Because $T$ is a regular Volterra operator with respect to $\xi,$ the operator $T|_{\xi_k(E)}$ is a regular Volterra operator with respect to $\xi'$. 

For the second proof, consider the collection $\mathcal{B}=\{R_t\colon t\in [0,1] \}$ appearing in the proof of Proposition~\ref{Transfer1}. Let us put $\mathcal{B}'\coloneqq \{R_t\cap (\xi_k(E)\times \xi_k(E))\colon t\in [0,1] \}$. It follows that domain the relation $R_t\cap (\xi_k(E)\times \xi_k(E))$ equals to the band $\xi_s(E)$ of $\xi_k(E)$ where $s=k$ if $\lfloor \phi^{-1}(t) \rfloor\geq k$ and $s=t$ otherwise. The operator $T|_{\xi_k(E)}\colon \xi_k(E)\to \xi_k(E)$ is a Volterra operator, in the sense of~\cite{ZA}, with respect to system $\mathcal{B}'$. Because $\phi\colon (0,\infty)\to (0,1)$ is order reversing, domains of $R_t\cap (\xi_k(E)\times \xi_k(E))$ form a nondecreasing system of bands of $\xi_k(E)$ as $t$ increases. This shows that $T$ is a regular Volterra with respect to some forward filtration in $\bb(\xi_k(E))$.
\end{proof}

\begin{prop}\label{PropVolttoRegular}
Suppose that $T\colon E\to E$ is a $\mathfrak{B}$-Volterra operator on a Banach lattice $E$ for some Boolean subalgebra $\mathfrak{B}$ of $\mathfrak{B}(E)$. Then $T$ is a regular Volterra operator with respect to every forward filtration in $\bb$. In particular, positive band preserving operators on $E$ are both regular Volterra, and, Volterra in the sense of~\cite{ZA}.
\end{prop}
\begin{proof}
	Result follows from definitions and from Corollary~\ref{CoroBandPreserving}.
\end{proof}

\begin{exam}
	In view of Example~\ref{Ex1} and Proposition~\ref{PropVolttoRegular}, a regular Volterra operator can be a compact operator. Indeed, in the settings of Proposition~\ref{PropVolttoRegular}, any compact $\bb$-Volterra operator is a compact regular Volterra operator.
\end{exam}

Recall from~\cite{KLW} that a positive order continuous projection $T$ on $E$ is said to be a conditional expectation if $Te$ is a weak order unit in $E$ whenever $e$ is a weak unit in $E$.
\begin{prop}\label{ConditionalExp}
	Let $T\colon E\to E$ be a conditional expectation, see~\cite{KLW}. Then $T$ is regular Volterra; and hence $T$ is Volterra in the sense of~\cite{ZA}. 
\end{prop} 
\begin{proof}
	Let $(x_n)_{n=1}^{\infty}$ be a sequence in $T(E)$ such that $\pi_{x_n}\leq \pi_{x_{n+1}}$ for all $n\geq 1$ where $\pi_{x_n}$ denotes the order projection $\pi_{x_n}(z)=\sup\{ z\wedge m|x_n|\colon m\in \nn \}$ for all $n\geq 1$. Such an increasing sequence of order projections exists because $Te\in T(E)$ is a weak order unit of $E,$ and hence $\pi_{Te}=\textbf{1}$ in $\bb(E)$. By Lemma 3.1 of~\cite{KLW}, the band $\pi_{x_n}(E)$ is left fixed by $T$ for all $n\geq 1$. This is equivalent to $\pi_{x_n}x=\pi_{x_n}y$ for $x,y\in E$ implies $\pi_{x_n}Tx=\pi_{x_n}Ty$ for all $n$. To see this, let $x-y\in \ker \pi_{x_n}$ so that $x-y\in \pi_{x_n}^*(E)$. By Lemma 3.1 of~\cite{KLW} again, $T(x-y)\in \pi_{x_n}^*(E)$. Hence, $\pi_{x_n} Tx=\pi_{x_n} Ty$. Hence, $T$ is a regular Volterra operator with respect to forward filtration \[\xi: \xi_0=\textbf{0}\leq \pi_{x_1}\leq\pi_{x_2}\leq \cdots\leq \xi_{\infty}=\textbf{1}.\] By Proposition~\ref{Transfer1}, T is a Volterra operator in the sense of~\cite{ZA}.
\end{proof}

\section{Actions on Regular Bounded Martingale Spaces}
\subsection{Induced Actions of $\bb$-Volterra Operators}\label{ForwardPaths}

Consider the function \[L\colon \mathcal{P}^+(\bb)\to \mathcal{P}^+(\bb)\] defined by $L(\xi)_0=\xi_0$ and $L(\xi)_n=\xi_{n+1}$ for $n\geq 1$. In tabular form, see Section~\ref{Paths}, if \[\xi\colon \xi_0=\textbf{0}\leq \xi_1\leq \xi_2\leq \cdots \leq \xi_{\infty}=\textbf{1}\]
then \[L^k(\xi)\colon \xi_0=\textbf{0}\leq \xi_{k+1}\leq \xi_{k+2}\leq \cdots \leq \xi_{\infty}=\textbf{1}\] for all $k\geq 0$. By convention, $L^0(\xi)=\xi$ for every $\xi \in \mathcal{P}^+(\bb)$.

In~\cite[Sec. 14]{T}, conditions on $T\colon E\to E$ with respect to which $T$ induces a linear action on the abstract martingale spaces are given.

\begin{lem}\label{P009}
Suppose that $T\colon E\to E$ is a $\mathfrak{B}$-Volterra operator on an order continuous Banach lattice $E$ for some Boolean subalgebra $\mathfrak{B}$ of $\mathfrak{B}(E)$. For every forward filtration $\xi\in \mathcal{P}^+(\mathfrak{B})$ the spaces $\mathcal{M}_0(\xi)$ and $\mathcal{M}_1(\xi^*)$ are isomorphic, and, the operator $T$ induces a positive operator $\hat{T}_{\xi}\colon \mathcal{M}_0(\xi)\to \mathcal{M}_0(\xi)$ making the diagrams
\begin{center}
	\begin{tikzcd}[ampersand replacement=\&]
	\mathcal{M}_0(\xi) \ar[r, "\hat{T}_{\xi}"]\ar[d, "\bold{s}"]\& \mathcal{M}_0(\xi)\ar[d, "\bold{s}"] \\ \mathcal{M}_0(L(\xi)) \ar[r, "\hat{T}_{L(\xi)}"]\& \mathcal{M}_0(L(\xi))
	\end{tikzcd}
	\quad \quad
	\begin{tikzcd}[ampersand replacement=\&]
	\mathcal{M}_1(\xi^*) \ar[r, "\hat{T}_{\xi}"]\ar[d, "\bold{s}"]\& \mathcal{M}_1(\xi^*)\ar[d, "\bold{s}"] \\ \mathcal{M}_1(L(\xi^*)) \ar[r, "\hat{T}_{L(\xi)}"]\& \mathcal{M}_1(L(\xi^*))
	\end{tikzcd}
\end{center}
commutative where $\textbf{s}\colon \mathcal{M}_i(\xi)\to \mathcal{M}_i(L(\xi))$ for $i=0,1$ denotes the forward shift operator defined by $(\textbf{s}((x_n)_{n=1}^\infty))_k=x_{k+1}$ for $k\geq 1$.
\end{lem}

\begin{proof}
Let $(x_n)_{n=1}^{\infty}\in \mathcal{M}_0(\xi)$ so that $\xi_nx_m=x_n$ holds for $m\geq n\geq 1$. Because $\xi$ is a forward filtration in $\bb,$ we have $\xi_n\in \bb$ for all $n$. Consider the operator $\hat{T}_{\xi}\colon \mathcal{M}_0(\xi)\to \mathcal{M}_0(\xi)$ defined by \[ \hat{T}_{\xi}((x_n)_{n=1}^{\infty})=(\xi_nTx_n)_{n=1}^{\infty}.\] We first observe that if $m\geq n$ then by Proposition~\ref{P001} \[\xi_nTx_m=\xi_nT\xi_m x_m=\xi_nT\xi_n\xi_m x_m=\xi_nT\xi_nx_m=\xi_nTx_n .\] Hence, \[\xi_n\xi_mTx_m=\xi_nTx_m=\xi_nTx_n. \] We conclude that $\hat{T}_{\xi}((x_n)_{n=1}^{\infty})\in \mathcal{M}_0(\xi)$. Evidently $\hat{T}_{\xi}$ is a positive operator.
	
Consider the operator $\textbf{s}\colon \mathcal{M}_0(\xi)\to \mathcal{M}_0(L(\xi))$ defined by \[(\textbf{s}((x_n)_{n=1}^\infty))_k=x_{k+1}\] for $k\geq 1$. If $x= (x_n)_{n=1}^\infty\in \mathcal{M}_0(\xi)$ then we have \[ L(\xi)_n(\textbf{s}(x)_m)=\xi_{n+1}x_{m+1}=x_{n+1}=\textbf{s}(x)_n\] for every $1\leq n\leq m$. Because \[(\textbf{s}\hat{T}_{\xi}x)_n=(\hat{T}_{\xi}x)_{n+1}=\xi_{n+1}Tx_{n+1}=L(\xi)_nT(\textbf{s}(x))_n=(\hat{T}_{L(\xi)}\textbf{s}(x))_n \] for every $n,$ the first diagram commutes.
	
It follows from definitions, see Section~\ref{Paths}, and from $\xi^{**}=\xi$ that the ordered vector spaces $\mathcal{M}_1(\xi^*)$ and $\mathcal{M}_0(\xi)$ are identical. Similarly, $\mathcal{M}_1(L(\xi^*))$ and $\mathcal{M}_0(L(\xi))$ are identical. These imply that the first diagram commutes if and only if the second diagram commutes.
\end{proof}

We see from Lemma~\ref{P009} that if $T$ is $\bb$-Volterra then $\hat{T}_{\xi}$ exists for every forward filtration $\xi$ in $\bb$. In the case of regular Volterra operators of Section~\ref{RegularVolterra}, $\hat{T}_{\xi}$ exists only for some forward filtration $\xi$. Details are provided below.

\begin{coro}
		Suppose that $T\colon E\to E$ is a regular Volterra operator with respect to some forward filtration $\xi$ in the Boolean algebra $\bb(E)$. Then $T$ induces an operator $\hat{T}_{\xi}\colon \mathcal{M}_0(\xi)\to \mathcal{M}_0(\xi)$ making the diagram
		\begin{center}
			\begin{tikzcd}[ampersand replacement=\&]
			\mathcal{M}_0(\xi) \ar[r, "\hat{T}_{\xi}"]\ar[d, "\bold{s}"]\& \mathcal{M}_0(\xi)\ar[d, "\bold{s}"] \\ \mathcal{M}_0(L(\xi)) \ar[r, "\hat{T}_{L(\xi)}"]\& \mathcal{M}_0(L(\xi))
			\end{tikzcd}
		\end{center}
		commutative. In particular, a conditional expectation $T\colon E\to E,$ see~\cite{KLW}, induces an action on $\mathcal{M}_0(\xi)$ for some forward filtration $\xi$.
\end{coro}
\begin{proof}
	Proof is similar to that of Lemma~\ref{P009}. Particular case follows from Proposition~\ref{ConditionalExp}. 
\end{proof}

\begin{prop}\label{PropGen}
	Suppose that $\xi\in \mathcal{P}^+(\mathfrak{B})$ is a forward filtration in a Boolean subalgebra $\bb$ of $\bb(E)$. The operator $T\mapsto \hat{T}_{\xi}$ sending a $\bb$-Volterra operator $T$ to the positive operator $\hat{T}_{\xi}$ is linear and order preserving.
\end{prop}
\begin{proof}
	Linearity is clear. Let $T,S$ be two $\bb$-Volterra operators such that $T\leq S$ with respect to the order, see~\cite[Chapter 1.1]{AB}, on the positive operators. If $(x_n)_{n=1}^{\infty}\in \mathcal{M}_0(\xi)$ is a positive abstract martingale, then \[((\hat{S}_{\xi}-\hat{T}_{\xi})(x_n)_{n=1}^{\infty})_k=\xi_k(S-T)(x_k)\geq 0\] for all $k\geq 1$. This implies that $\hat{T}_{\xi}\leq \hat{S}_{\xi}$ on positive abstract martingales.
\end{proof}
\begin{coro}\label{Lemma2}
	Suppose that $\xi\in \mathcal{P}^+(\mathfrak{B})$ is a forward filtration in a Boolean subalgebra $\bb$ of $\bb(E)$. 
\begin{itemize}
		\item[\em i.]{If $T,S\colon E\to E$ are commuting $\bb$-Volterra operators then the operators $\hat{T}_{\xi}$ and $\hat{S}_{\xi}$ on $\mathcal{M}_0(\xi)$ also commute. }
		\item[\em ii.]{If $T$ is $\bb$-Volterra and $S=\pi T$ for some $\pi\in \bb$ then $\hat{S}_{\xi}=\pi\cdot \hat{T}_{\xi}$.}
		\item[\em iii.]{If $T$ is $\bb$-Volterra and $S=\pi+T-\pi T$ for some $\pi\in\bb$ then $\hat{S}_{\xi}=\hat{\pi}_{\xi}+\hat{T}_{\xi}-\pi\cdot\hat{T}_{\xi}$.}
\end{itemize}	
\end{coro}
\begin{proof}
$(i)$. Let $(x_n)_{n=1}^{\infty}\in \mathcal{M}_0(\xi)$ be an abstract martingale. Then
\begin{align*}
		\hat{T}_{\xi}(\hat{S}_{\xi}((x_n)_{n=1}^{\infty}))&=(\xi_n T\xi_nSx_n)_{n=1}^{\infty}=(\xi_nTSx_n)_{n=1}^{\infty}\\
		&=(\xi_n S\xi_nTx_n)_{n=1}^{\infty}=\hat{S}_{\xi}(\hat{T}_{\xi}((x_n)_{n=1}^{\infty}))
\end{align*}
by Proposition~\ref{P001} because both $T$ and $S$ are $\bb$-Volterra.

$(ii)$. Because $T$ is $\bb$-Volterra, for every $\pi\in \bb$ the operator $\pi T$ is $\bb$-Volterra by Lemma~\ref{Lemma1}. As $\bb$ is a Boolean algebra, order projections belonging to $\bb$ commute. Hence, \[ \hat{S}_{\xi}((x_n)_{n=1}^{\infty})=(\xi_n \pi Tx_n)_{n=1}^{\infty}=\pi(\xi_n Tx_{n})_{n=1}^{\infty}=\pi\cdot \hat{T}_{\xi}((x_n)_{n=1}^{\infty})\] for $(x_n)_{n=1}^{\infty}\in \mathcal{M}_0(\xi)$. 

$(iii)$. Because $T$ is $\bb$-Volterra, the operator $\pi+T-\pi T$ is $\bb$-Volterra for every $\pi\in \bb$ by Lemma~\ref{Lemma1}. Hence, \[ \hat{S}_{\xi}((x_n)_{n=1}^{\infty})=(\pi\xi_nx_n+\xi_nTx_n-\xi_n\pi Tx_n)_{n=1}^{\infty}=(\hat{\pi}_{\xi}+\hat{T}_{\xi}-\pi\cdot\hat{T}_{\xi})((x_n)_{n=1}^{\infty})\] for $(x_n)_{n=1}^{\infty}\in \mathcal{M}_0(\xi)$. 
\end{proof}

\begin{coro} Suppose that $\xi\in \mathcal{P}^+(\mathfrak{B})$ is a forward filtration in a Boolean subalgebra $\bb$ of $\bb(E)$. Let $T$ be a $\bb$-Volterra operator.
\begin{itemize}
		\item[\em i.]{If $T\leq I_E$ and $S=\pi-\pi T+\pi^* T$ for some $\pi\in B$ then $\hat{S}_{\xi}=\hat{\pi}_{\xi}-\pi\cdot \hat{T}_{\xi}+\pi^* \cdot\hat{T}_{\xi}$ is a positive operator on $\mathcal{M}_0(\xi)$.}
		\item[\em ii.]{If $\pi\in \bb$ satisfies $\pi T=\pi$ and $S=I_E+\pi-T$ then $\hat{S}_{\xi}=\hat{(I_E)}_{\xi}+\hat{\pi}_{\xi}-\hat{T}_{\xi}$ is a positive operator on $\mathcal{M}_0(\xi)$.}
\end{itemize}	
\end{coro}
\begin{proof}
$(i)$. If $T\leq I_E$ then by Theorem~\ref{P0010}, the operator $S=\pi-\pi T+\pi^* T$ is $\bb$-Volterra, in particular it is positive. By Proposition~\ref{PropGen} and Corollary~\ref{Lemma2}, the operator $\hat{S}_{\xi}=\hat{\pi}_{\xi}-\pi \hat{T}_{\xi}+\pi^* \hat{T}_{\xi}$ is positive. 

$(ii)$. The proof is similar to the proof of $(i),$ and it follows from Theorem~\ref{P0010}, Corollary~\ref{Lemma2} and Proposition~\ref{PropGen}.
\end{proof}

We first note from~\cite[Example 2.4]{SG1990} and ~\cite[Example 1]{SG2019} that for each vector in $E$ there corresponds an abstract martingale. In the case of bounded martingales, by~\cite[Sec.9]{T} and~\cite[Sec.3]{GT}, for every forward filtration $\xi$ there is a contraction $\iota_{\xi}\colon E\to \mathcal{M}_{b}(\xi)$ defined by $\iota_{\xi}(x)=(\xi_n x)_{n=1}^{\infty}$. 
 
\begin{thm}\label{P0KB}
 	Let $E$ be a $KB$-space. Suppose that $\xi\in \mathcal{P}^+(\mathfrak{B})$ is a forward filtration in a Boolean subalgebra $\bb$ of $\bb(E)$. Let $T$ be a $\bb$-Volterra operator. The operator $\hat{T}_{\xi}\colon \mathcal{M}_0(\xi)\to \mathcal{M}_0(\xi)$ induces a norm continuous operator $\hat{T}_{\xi}\colon \mathcal{M}_{b}(\xi)\to \mathcal{M}_{b}(\xi)$. Denote by $\iota_{\xi}\colon E\to \mathcal{M}_{b}(\xi)$ the contraction given in~\cite[Sec.9]{T}. Then the diagram 
 	\begin{center}
 		\begin{tikzcd}[ampersand replacement=\&]
 		E\ar[r, "\iota_{\xi}"]\ar[d, "T"]\& \mathcal{M}_{b}(\xi)\ar[d, "\hat{T}_{\xi}"]\ar[r, "\bold{s}"] \& \mathcal{M}_{b}(L(\xi))\ar[d, "\hat{T}_{L(\xi)}"]\\ E \ar[r, "\iota_{\xi}"]\& \mathcal{M}_{b}(\xi)\ar[r, "\bold{s}"]\&\mathcal{M}_{b}(L(\xi))
 		\end{tikzcd}
 	\end{center}
of norm continuous operators and Banach lattices commutes.
 \end{thm}
\begin{proof}
By Lemma~\ref{P009}, the operator $\hat{T}_{\xi}\colon \mathcal{M}_0(\xi)\to \mathcal{M}_0(\xi)$ is positive. By Theorem 1.10 of~\cite{AB}, it induces a positive operator $\hat{T}_{\xi}\colon \mathcal{M}_{b}(\xi)\to \mathcal{M}_{b}(\xi)$. As $E$ is a KB-space, by Theorem 7 of~\cite{T}, the space $\mathcal{M}_{b}(\xi)$ is a Banach lattice. By Theorem 4.3 of~\cite{AB}, $\hat{T}_{\xi}\colon \mathcal{M}_{b}(\xi)\to \mathcal{M}_{b}(\xi)$ is continuous with respect to norm topology induced by the canonical norm on the space of bounded martingales.

In order to show commutativity of the diagram, recall from~\cite[Sec.9]{T} that $\iota_{\xi}\colon E\to \mathcal{M}_{b}(\xi)$ is given by $\iota_{\xi}(x)=(\xi_n x)_{n=1}^{\infty}$. It follows from Proposition~\ref{P001} that \[(\iota_{\xi}\circ T)(x)=(\xi_n Tx)_{n=1}^{\infty}=(\xi_nT\xi_nx)_{n=1}^{\infty}=\hat{T}_{\xi}(\xi_nx)_{n=1}^{\infty}=(\hat{T}_{\xi}\circ \iota_{\xi})(x) \] for every $x\in E$. The second square of the diagram commutes by Lemma~\ref{P009}.
\end{proof}

\begin{coro} In the settings of Theorem~\ref{P0KB}, the diagrams 
		\begin{center}
		\begin{tikzcd}[ampersand replacement=\&]
		E\ar[r, "\iota_{L^k(\xi)}"]\ar[d, "T"]\& \mathcal{M}_{b}(L^k(\xi))\ar[d, "\hat{T}_{L^k(\xi)}"]\ar[r, "\bold{s}"] \& \mathcal{M}_{b}(L^{k+1}(\xi))\ar[d, "\hat{T}_{L^{k+1}(\xi)}"]\\ E \ar[r, "\iota_{L^k(\xi)}"]\& \mathcal{M}_{b}(L^k(\xi))\ar[r, "\bold{s}"]\&\mathcal{M}_{b}(L^{k+1}(\xi))
		\end{tikzcd}
		\end{center} and
		\begin{center}
		\begin{tikzcd}[ampersand replacement=\&]
		E\ar[r, "\iota_{L^k(\xi)}"] \ar[dr, "\iota_{L^{k+1}(\xi)}"']\& \mathcal{M}_{b}(L^k(\xi))\ar[d, "\bold{s}"]\\ \& \mathcal{M}_{b}(L^{k+1}(\xi))
		\end{tikzcd}
	\end{center}
of norm continuous operators and Banach lattices commute for all $k\geq 0$. Moreover, if $\hat{T}_{\xi}$ is surjective then $\hat{T}_{L^{k}(\xi)}$ is open for every $k\geq 1$.
\end{coro}
\begin{proof}
	Commutativity of the first diagram is a direct corollary of Theorem~\ref{P0KB} because $L^k(\xi)$ is a forward filtration for $k\geq 0$ whenever $\xi$ is a forward filtration. In order to verify commutativity of the second triangular diagram, let $x\in E$ be arbitrary. We have \[ \bold{s}(\iota_{L^k(\xi)}(x))=\bold{s}(\xi_{k+n}x)_{n=1}^{\infty}=(\xi_{k+n+1}x)_{n=1}^{\infty}=\iota_{L^{k+1}(\xi)}(x)\] for all $k\geq 0$.
	
	Because the operators $\bold{s}$ and $\hat{T}_{\xi}$ are surjective, they are open by the open mapping theorem. Composition of open maps is open. Commutativity of the first diagram implies that $\hat{T}_{L^{k}(\xi)}$ is open for every $k\geq 1$.
\end{proof}

\begin{exam} Let $\pi\colon L_p[0,1]\to L_p[0,1]$ be $\pi(f)=f\cdot \chi_{[0,1/2]}$ where $\chi_{[0,1/2]}$ denotes the characteristic function of $[0,1/2]$. Let $\bb=\{ \textbf{0},\textbf{1},\pi,\pi^*\}$. Consider the forward filtration given by the tabular form \[\xi\colon \xi_0=\textbf{0}\leq \xi_1=\pi\leq \xi_2=\pi\leq \cdots\leq \xi_{\infty}=\textbf{1}. \] We have $(f_n)_{n=1}^\infty\in \mathcal{M}_{b}(\xi)$ if $f_1=\pi(f_n)=f_n$ for all $n\geq 1$. Thus $(f_n)_{n=1}^\infty\in \mathcal{M}_{b}(\xi)$ if and only if it is of the form $(f_n)_{n=1}^\infty=f_1,f_1,f_1,\ldots$ where $f_1$ is equivalent to zero outside $[0,1/2]$. The operator $\pi/2$ is $\bb$-Volterra. The induced operator $\hat{(\frac{\pi}{2})}_{\xi}\colon \mathcal{M}_{b}(\xi)\to \mathcal{M}_{b}(\xi)$ of Theorem~\ref{P0KB} is given by $\hat{(\frac{\pi}{2})}_{\xi}(x)=x/2$ for $x\in \mathcal{M}_{b}(\xi)$. 
\end{exam}

We recall the convention that $\textbf{s}^0$ and $L^0$ are equal to the identity maps on their domains of definitions. 

\begin{thm} Let $E$ be a $KB$-space. Suppose that $\xi\in \mathcal{P}^+(\mathfrak{B})$ is a forward filtration in a Boolean subalgebra $\bb$ of $\bb(E)$. Let $T\colon E\to E$ be a $\bb$-Volterra operator. Denote by $B'$ the band generated by \[ \{(\textbf{s}^k(x))_{k=0}^{\infty}\colon x\in \mathcal{M}_{b}(\xi) \} \subseteq (\bigoplus_{k=0}^{\infty}\mathcal{M}_{b}(L^k(\xi)))_{\ell_{\infty}}\] in the Banach lattice $(\bigoplus_{k=0}^{\infty}\mathcal{M}_{b}(L^k(\xi)))_{\ell_{\infty}}$. For each $l\geq 1,$ the operator \[(\hat{T}_{\xi}^l)^{!}((y_k)_{k=0}^{\infty})=(\textbf{s}^{k}\hat{T}_{\xi}^ly_0)_{k=0}^{\infty} \] for $(y_k)_{k=0}^{\infty}\in (\bigoplus_{k=0}^{\infty}\mathcal{M}_{b}(L^k(\xi)))_{\ell_{\infty}}$ extends to a norm continuous operator \[(\hat{T}_{\xi}^l)^{!}\colon B'\to (\bigoplus_{k=0}^{\infty}\mathcal{M}_{b}(L^k(\xi)))_{\ell_{\infty}}\] such that $(\hat{T}_{\xi}^l)^{!}=((\hat{T}_{\xi}^1)^{!})^l,$ i.e., the correspondence $!$ is a discrete semigroup morphism between the semigroups generated by $\hat{T}_{\xi}$ and $(\hat{T}_{\xi}^1)^{!}$. 
\end{thm}
\begin{proof} As $E$ is a $KB$-space, $(\bigoplus_{k=0}^{\infty}\mathcal{M}_{b}(L^k(\xi)))_{\ell_{\infty}}$ is a Banach lattice under pointwise ordering by~\cite[Sec.4.1]{AB}. The set $\{(\textbf{s}^k(x))_{k=0}^{\infty}\colon x\in \mathcal{M}_{b}(\xi) \}$ is a subset of $(\bigoplus_{k=0}^{\infty}\mathcal{M}_{b}(L^k(\xi)))_{\ell_{\infty}}$ because \[\sup_{k\geq 0}||\textbf{s}^k(x)||_{\mathcal{M}_{b}(L^k(\xi))}\leq ||x||_{\mathcal{M}_{b}(\xi)}<\infty \] for every $x\in \mathcal{M}_{b}(\xi)$. Thus, the set $\{(\textbf{s}^k(x))_{k=0}^{\infty}\colon x\in \mathcal{M}_{b}(\xi) \}$ generates a band in $(\bigoplus_{k=0}^{\infty}\mathcal{M}_{b}(L^k(\xi)))_{\ell_{\infty}}$. Denote this band by $B'$. It follows from the commutative diagram given in Theorem~\ref{P0KB} that
\begin{align*}
(\hat{T}_{\xi}^l)^{!}((y_k+z_k)_{k=0}^{\infty})&=(\textbf{s}^{k}\hat{T}_{\xi}^l(y_0+z_0))_{k=0}^{\infty} \\ &= (\hat{T}_{\xi}^ly_k)_{k=0}^{\infty}+ (\hat{T}_{\xi}^lz_k)_{k=0}^{\infty} \\&= (\textbf{s}^{k} \hat{T}_{\xi}^ly_0)_{k=0}^{\infty}+ (\textbf{s}^{k} \hat{T}_{\xi}^lz_0)_{k=0}^{\infty} \\ &=(\hat{T}_{\xi}^l)^{!}(y_k)_{k=0}^{\infty}+(\hat{T}_{\xi}^l)^{!}(z_k)_{k=0}^{\infty}
\end{align*}	
for $(y_k)_{k=0}^{\infty},(z_k)_{k=0}^{\infty}\in (B')^+$ and $l\geq 0$. This implies that $(\hat{T}_{\xi}^k)^{!}$ is additive on $(B')^+$. Hence, by Theorem 1.10 of~\cite{AB} we have the result.

Because $T$ is $\bb$-Volterra $(\hat{T}_{\xi})^l=\hat{(T^l)_{\xi}}$ for every $l\geq 0$. This implies that the correspondence~$!$ is a semigroup morphism between the semigroups generated by $\hat{T}_{\xi}$ and $(\hat{T}_{\xi}^1)^{!}$. 
\end{proof}

\begin{thm}\label{Fundcoro}
Suppose that $T\colon E\to E$ is a $\mathfrak{B}$-Volterra operator on an order continuous Banach lattice $E$ for some Boolean subalgebra $\mathfrak{B}$ of $\mathfrak{B}(E)$. The operator $\hat{T}_{\xi}\colon \mathcal{M}_0(\xi)\to \mathcal{M}_0(\xi)$ induces a norm continuous operator $\hat{T}_{\xi}\colon \mathcal{M}_{r}(\xi)\to \mathcal{M}_{r}(\xi),$ again denoted by the same notation, on the Banach lattice $\mathcal{M}_{r}(\xi)$ of all regular bounded martingales with respect to regular norm, see Section~\ref{AbtractMartingaleSpaces}. The diagram 
\begin{center}
	\begin{tikzcd}[ampersand replacement=\&]
	\mathcal{M}_{r}(\xi) \ar[r, "\hat{T}_{\xi}"]\ar[d, "\bold{s}"]\& \mathcal{M}_{r}(\xi)\ar[d, "\bold{s}"] \\ \mathcal{M}_{r}(L(\xi)) \ar[r, "\hat{T}_{L(\xi)}"]\& \mathcal{M}_{r}(L(\xi))
	\end{tikzcd}
\end{center}
of Banach lattices and norm continuous operators is commutative up to an equivalence of Banach lattices.
\end{thm}
\begin{proof}
	The operator $\hat{T}_{\xi}\colon \mathcal{M}_0(\xi)\to \mathcal{M}_0(\xi)$ is positive. By Theorem 1.10 of~\cite{AB}, it induces a positive operator $\hat{T}_{\xi}\colon \mathcal{M}_{r}(\xi)\to \mathcal{M}_{r}(\xi)$. By Theorem 3.5 of~\cite{TX}, the space $\mathcal{M}_{r}(\xi)$ is a Banach lattice under regular norm. By Theorem 4.3 of~\cite{AB}, $\hat{T}_{\xi}\colon \mathcal{M}_{r}(\xi)\to \mathcal{M}_{r}(\xi)$ is continuous with respect to regular norm. Because the diagram commutes, it also commutes up to an equivalence of Banach lattices.
\end{proof}

In Theorem~\ref{Fundcoro} we obtain a commutative diagram of Banach lattices and norm continuous operators. New Banach spaces and operators can be deduced via this commutative diagram. Details are given in Definition~\ref{DefAssociated}.

\begin{definition}\label{DefAssociated}
		Suppose that $T\colon E\to E$ is a $\mathfrak{B}$-Volterra operator on an order continuous Banach lattice $E$ for some Boolean subalgebra $\mathfrak{B}$ of $\mathfrak{B}(E)$. Let $\xi$ be a forward filtration in $\bb$. A Banach space $\mathcal{M}_{T,\xi}$ is said to be an associated space of the pair $(T,\xi)$ if $\mathcal{M}_{T,\xi}$ is isometrically Banach space isomorphic to a direct limit of the directed system
		\begin{center}
		\begin{tikzcd}[ampersand replacement=\&]
		\mathcal{M}_{r}(\xi) \ar[r, "\hat{T}_{\xi}"] \& \mathcal{M}_{r}(\xi) \ar[r, "\hat{T}_{\xi}"] \& \mathcal{M}_{r}(\xi) \ar[r, "\hat{T}_{\xi}"] \& \cdots.
		\end{tikzcd}
		\end{center}
		 In this case, we write \[\mathcal{M}_{T,\xi}\coloneqq \varinjlim (\mathcal{M}_{r}(\xi), \hat{T}_{\xi} ).\] 
		 
		 A Banach space $\mathcal{M}_{\xi}$ is said to be an associated space of the forward filtration $\xi$ if $\mathcal{M}_{\xi}$ is isometrically Banach space isomorphic to a direct limit of the directed system 
		\begin{center}
			\begin{tikzcd}[ampersand replacement=\&]
			\mathcal{M}_{r}(\xi) \ar[r, "\bold{s}"] \& \mathcal{M}_{r}(L(\xi)) \ar[r, "\bold{s}"] \& \mathcal{M}_{r}(L^2(\xi)) \ar[r, "\bold{s}"] \& \cdots
			\end{tikzcd}
		\end{center}
	 In this case, we write $\mathcal{M}_{\xi}\coloneqq \varinjlim_k (\mathcal{M}_{r}(L^k(\xi)), \bold{s} )$.
\end{definition}

As a result of Theorem~\ref{Fundcoro} we obtain the following. 

\begin{coro}\label{AssociatedProp}
	Suppose that $T\colon E\to E$ is a $\mathfrak{B}$-Volterra operator on an order continuous Banach lattice $E$ for some Boolean subalgebra $\mathfrak{B}$ of $\mathfrak{B}(E)$. Let $\xi$ be a forward filtration in $\bb$. 
	\begin{itemize}
		\item[\em i.]{ The forward shift map $\mathbf{s}\colon \mathcal{M}_{r}(\xi)\to \mathcal{M}_{r}(L(\xi))$ induces a norm continuous operator \[ \hat{\bold{s}^k}\colon \mathcal{M}_{T,L^{k-1}(\xi)}\to \mathcal{M}_{T,L^{k}(\xi)}\] from the associated space $\mathcal{M}_{T,L^{k-1}(\xi)}$ of the pair $(T, L^{k-1}(\xi))$ to the associated space $\mathcal{M}_{T,L^{k}(\xi)}$ of the pair $(T,L^{k}(\xi))$ for every $k\geq 1$. }
		\item[\em ii.]{The operator $T\colon E\to E$ induces a norm continuous operator $\hat{\hat{T}}\colon \mathcal{M}_{\xi}\to \mathcal{M}_{\xi}$ on the associated space $\mathcal{M}_{\xi}$ of $\xi$. }
	\end{itemize}
\end{coro}

\begin{proof}
$(i)$. Because the diagram given in Theorem~\ref{Fundcoro} is commutative, the diagram 
\begin{center}
	\begin{tikzcd}[ampersand replacement=\&]
	\mathcal{M}_{r}(L^{k-1} (\xi)) \ar[r, "\hat{T}_{L^{k-1}(\xi)}"]\ar[d, "\bold{s^k}"]\& \mathcal{M}_{r}(L^{k-1} (\xi))\ar[d, "\bold{s^k}"] \\ \mathcal{M}_{r}(L^{k} (\xi)) \ar[r, "\hat{T}_{L^k(\xi)}"]\& \mathcal{M}_{r}(L^{k} (\xi))
	\end{tikzcd}
\end{center}
commutes for all $k\geq 1$. We remark that $\bold{s^k}$ denotes the forward shift map $\bold{s}$ with domain $\mathcal{M}_{r}(L^{k-1} (\xi))$. By considering definitions of $\mathcal{M}_{T,L^{k-1}(\xi)}$ and $\mathcal{M}_{T,L^{k}(\xi)}$, we have
\begin{center}
		\begin{tikzcd}[ampersand replacement=\&]
		\mathcal{M}_{r}(L^{k-1} (\xi))\ar[r, "\hat{T}_{L^{k-1}(\xi)}"]\ar[d, "\bold{s^k}"]\& \mathcal{M}_{r}(L^{k-1} (\xi))\ar[d, "\bold{s^k}"]\ar[r, "\hat{T}_{L^{k-1}(\xi)}"] \& \mathcal{M}_{r}(L^{k-1} (\xi))\ar[d, "\bold{s^k}"]\ar[r,"\hat{T}_{L^{k-1}(\xi)}"] \&\cdots \\ \mathcal{M}_{r}(L^{k} (\xi)) \ar[r, "\hat{T}_{L^k(\xi)}"]\& \mathcal{M}_{r}(L^{k} (\xi))\ar[r, "\hat{T}_{L^k(\xi)}"]\&\mathcal{M}_{r}(L^{k} (\xi))\ar[r, "\hat{T}_{L^k(\xi)}"] \& \cdots
		\end{tikzcd}
		\end{center} from which the result follows.
		
$(ii)$. The diagram
\begin{center} 
\begin{tikzcd}[ampersand replacement=\&]
		\mathcal{M}_{r}(\xi)\ar[r, "\bold{s}"]\ar[d, "\hat{T}_{\xi}"]\& \mathcal{M}_{r}(L (\xi))\ar[d, "\hat{T}_{L(\xi)}"]\ar[r, "\bold{s}"] \& \mathcal{M}_{r}(L^{2} (\xi))\ar[d, "\hat{T}_{L^2(\xi)}"]\ar[r,"\bold{s}"] \&\cdots \\ \mathcal{M}_{r}(\xi) \ar[r, "\bold{s}"]\& \mathcal{M}_{r}(L(\xi))\ar[r, "\bold{s}"]\&\mathcal{M}_{r}(L^{2} \xi)\ar[r, "\bold{s}"] \& \cdots
		\end{tikzcd}
		\end{center}
commutes by Theorem~\ref{Fundcoro}.
\end{proof}

Our next result is based on sectionally open transformations.

\begin{thm}\label{StarOpenInduced}
	Let $\bb$ be a Boolean subalgebra of $\bb(E)$ where $E$ is an order continuous Banach lattice. Suppose that $S\colon \bb\to \bb$ is a sectionally open transformation fixing $\textbf{0}\in \bb,$ see Proposition~\ref{StarOpen}. Then $S$ induces a norm continuous operator $\hat{S}\colon \mathcal{M}_{S(\xi)}\to \mathcal{M}_{\xi}$ between the associated spaces of the forward filtrations $S(\xi)$ and $\xi$, respectively.
\end{thm}
\begin{proof}
Because $\xi$ is a forward filtration, $S(\xi)$ is a forward filtration by Proposition~\ref{StarOpen}. Recall from Definition~\ref{DefAssociated} that \[\mathcal{M}_{\xi}=\varinjlim_k (\mathcal{M}_{r}(L^k(\xi)), \bold{s})\] and that \[\mathcal{M}_{S(\xi)}=\varinjlim_k (\mathcal{M}_{r}(L^k(S(\xi))), \bold{s}).\] 
For the forward filtration $\xi\colon \xi_0=\textbf{0}\leq \xi_1\leq\cdots\leq \xi_{\infty}=\textbf{1}$ we have \[L^k(S(\xi))\colon S(\xi)_0=\textbf{0}\leq S(\xi_{k+1})\leq S(\xi_{k+2}) \cdots \leq S(\xi_{\infty})=\textbf{1} \]	and \[L^k(\xi)\colon \xi_0=\textbf{0}\leq \xi_{k+1}\leq \xi_{k+2} \cdots \leq \xi_{\infty}=\textbf{1} \] for all $k\geq 1$.

We claim that if $x\in \mathcal{M}_{r}(L^k(S(\xi)))$ is a regular bounded martingale, then we have $x\in \mathcal{M}_{r}(L^k(\xi))$. In order to verify the claim, we write $x$ as a difference of two positive bounded martingales, and verify the same claim for each of these positive bounded martingales. Let $x^1,x^2\in (\mathcal{M}_{b}(L^k(S(\xi))))^+$ be such that $x=x^1-x^2$ where $x^i=(x_n^i)_{n=1}^{\infty}$ for $i=1,2$. Because $x^i\in (\mathcal{M}_{b}(L^k(S(\xi))))^+$, we have \[S(\xi_{k+n})(x_{n+m}^i)=x_n^{i} \] for all $m\geq 0,$ $n\geq 1$ and $i=1,2$. Because $S$ is sectionally open, $S(\xi_{k+n})\geq \xi_{k+n}$. Hence, we have \[\xi_{k+n}(x_{n+m}^i)=x_n^{i} \] for all $m\geq 0,$ $n\geq 1$ and $i=1,2$. Hence, $x^1,x^2\in (\mathcal{M}_{b}(L^k(\xi)))^+$ and $x\in \mathcal{M}_{r}(L^k(\xi))$. Denote by $\iota_k\colon \mathcal{M}_{r}(L^k(S(\xi)))\to \mathcal{M}_{r}(L^k(\xi))$ the identity operator. The diagram 
\begin{center}
	\begin{tikzcd}[ampersand replacement=\&]
	\mathcal{M}_{r}(L^k(S(\xi))) \ar[r, "\bold{s}"]\ar[d, "\iota_k"]\& \mathcal{M}_{r}(L^{k+1}(S(\xi)))\ar[d, "\iota_{k+1}"] \\ \mathcal{M}_{r}(L^k(\xi)) \ar[r, "\bold{s}"]\& \mathcal{M}_{r}(L{k+1}(\xi))
	\end{tikzcd}
\end{center}
is commutative. Hence, the maps $(\iota_k)$ induces a map $\hat{S}\colon \mathcal{M}_{S(\xi)}\to \mathcal{M}_{\xi}$.
\end{proof}

\subsection{An application of antichains in $\bb$.}\label{Antichains} We recall from Section~\ref{Intro} that a nonempty subset $S\subseteq \bb$ is said to be an \textit{antichain} in $\bb$ if $\pi,\rho\in S$ are distinct elements then $\pi\wedge\rho=0$. 

\begin{definition}
We say that a nonempty set $\{\xi^1,\xi^2,\ldots\}\subseteq \mathcal{P}^+(\mathfrak{B})$ forms an \textit{antichain of forward filtrations} if the set $\{ \xi_n^1,\xi_n^2,\ldots,\}\subseteq \bb$ of order projections is an antichain in $\bb$ for every $1\leq n<\infty$. 
\end{definition}

Let $\xi^1$ and $\xi^2$ be two forward filtrations. A norm continuous operator $S\colon \mathcal{M}_{r}(\xi^1) \to \mathcal{M}_{r}(\xi^2)$ intertwines $\hat{T}_{\xi^1}$ on $\mathcal{M}_{r}(\xi^1)$ with $\hat{T}_{\xi^2}$ on $\mathcal{M}_{r}(\xi^2)$ if \[S\circ \hat{T}_{\xi^1}=\hat{T}_{\xi^2}\circ S.\] We denote by \[\mathcal{L}_T(\xi^1,\xi^2)\subseteq \mathcal{L}(\mathcal{M}_{r}(\xi^1),\mathcal{M}_{r}(\xi^2))\] the set of all norm continuous operators $S\colon \mathcal{M}_{r}(\xi^1) \to \mathcal{M}_{r}(\xi^2)$ intertwining $\hat{T}_{\xi^1}$ with $\hat{T}_{\xi^2}$.

\begin{coro}
	Suppose that $T\colon E\to E$ is a $\mathfrak{B}$-Volterra operator on an order continuous Banach lattice $E$ for some Boolean subalgebra $\mathfrak{B}$ of $\mathfrak{B}(E)$. For every forward filtration $\xi$ in $\bb$, the forward shift operator $\textbf{s}\colon \mathcal{M}_{r}(\xi) \to \mathcal{M}_{r}(L(\xi))$ belongs to $\mathcal{L}_T(\xi,L(\xi))$.
\end{coro}
\begin{proof}
	This is a direct corollary of Theorem~\ref{Fundcoro}.
\end{proof}
	
\begin{prop}
	Suppose that $T\colon E\to E$ is a $\mathfrak{B}$-Volterra operator on a Banach lattice $E$ for some Boolean subalgebra $\mathfrak{B}$ of $\mathfrak{B}(E)$. Let $\{ \xi^1,\xi^2\}\subseteq \mathcal{P}^+(\mathfrak{B})$ be an antichain of forward filtrations. Let $S_k\colon E\to E$ be a positive operator such that one of 
	 \begin{itemize}
	 	\item[\em i.] {$S_kT=\xi_k^2TS_k+S_k(\xi_k^1)^*T$}
	 	\item[\em ii.]{$TS_k=(\xi_k^2)^*TS_k+S_k\xi_k^1T$}
	 \end{itemize}
holds for each $k\geq 1$. Then the positive operator $S=\prod_{k=1}^\infty S_k$ belongs to $\mathcal{L}_T(\xi^1,\xi^2)$. 
\end{prop}
\begin{proof}
	We first show that $(i)$ if and only if $(ii)$ so that if one of these identities hold for some $k$ then the other also holds. If $(i)$ then $\xi_k^2TS_k=S_k(1-(\xi_k^1)^*)T=S_k\xi_k^1T$. If $(ii)$ then $S_k\xi_k^1T=(1-(\xi_k^2)^*)TS_k=\xi_k^2TS_k$. Therefore, $(i)$ if and only if $(ii)$ for each $k\geq 1$.
	Because $S_k$ is positive for all $k\geq 1,$ the product $S=S_1\times S_2\times\cdots $ is a positive operator on the Banach lattice $\mathcal{M}_{r}(\xi^1)$. In particular, it is continuous with respect to regular norm. The identities $(i)$ and $(ii)$ imply that $S\in \mathcal{L}_T(\xi^1,\xi^2)$. Indeed, \[ \hat{T}_{\xi^2}\circ S((x_n)_{n=1}^{\infty})=(\xi_n^2TS_nx_n)_{n=1}^\infty=(S_n\xi_n^1 Tx_n)_{n=1}^{\infty}= S\circ \hat{T}_{\xi^1}((x_n)_{n=1}^{\infty}) \] for every positive regular bounded martingale $(x_n)_{n=1}^{\infty}$.
\end{proof}

Recall from Section~\ref{ForwardPaths} that $L\colon \mathcal{P}^+(\bb)\to \mathcal{P}^+(\bb)$ is defined by $L(\xi)_0=\xi_0$ and $L(\xi)_n=\xi_{n+1}$ for $n\geq 1$.
\begin{prop}
	If $S_1\in \mathcal{L}_T(\xi^1,\xi^2)$ and $S_2\in \mathcal{L}_T(\xi^2,\xi^3)$ for $\xi^1,\xi^2,\xi^3\in \mathcal{P}^+(\mathfrak{B})$ then $S_2S_1\in \mathcal{L}_T(\xi^1,\xi^3)$. For every $S\in \mathcal{L}_T(\xi^1,\xi^2)$ there is a norm continuous operator $\hat{L}(S)\in \mathcal{L}_T(L(\xi^1),L(\xi^2))$ such that $\hat{L}(S)\circ \bold{s}=\bold{s}\circ S$. 
\end{prop}
\begin{proof}
	If $S_1\in \mathcal{L}(\xi^1,\xi^2)$ and $S_2\in \mathcal{L}(\xi^2,\xi^3)$ then $S_1 \hat{T}_{\xi^1}=\hat{T}_{\xi^2} S_1$ and $S_2\hat{T}_{\xi^2}=\hat{T}_{\xi^3}S_2$. Hence, $S_2S_1\hat{T}_{\xi^1}=S_2\hat{T}_{\xi^2}S_1=\hat{T}_{\xi^3}S_2S_1$. 
	
	For the second part, consider $S\in \mathcal{L}_T(\xi^1,\xi^2),$ together with $(x_n)_{n=1}^{\infty} \in \mathcal{M}_{r}(\xi^1)$ and $(y_n)_{n=1}^{\infty} \in \mathcal{M}_{r}(\xi^2)$ such that $S(x_n)_{n=1}^{\infty}=(y_n)_{n=1}^{\infty}$. We claim that if we erase the first coordinates of $(x_n)_{n=1}^{\infty}$ and $(y_n)_{n=1}^{\infty}$ then we obtain regular bounded martingales, and moreover, this positive operation induces a continuous operator satisfying the required properties.
	
	Formally, let us put $z_n=x_{n+1}$ and $w_n=y_{n+1}$ for $n\geq 1$. Evidently $(z_n)_{n=1}^{\infty}$ and $(w_n)_{n=1}^{\infty}$ belong to $\mathcal{M}_{r}(L(\xi^1))$ and $\mathcal{M}_{r}(L(\xi^2)),$ respectively. Conversely, for an arbitrary element $(z_n)_{n=1}^{\infty}$ of $\mathcal{M}_{r}(L(\xi^1))$ if we put $x_1=\xi_1^1z_1,$ $x_n=z_{n-1}$ for $n\geq 2$ then we obtain $(x_n)_{n=1}^{\infty}\in \mathcal{M}_{r}(\xi)$. 
	
	Hence given $S\in \mathcal{L}_T(\xi^1,\xi^2)$ we define $\hat{L}(S)$ by \[ S((x_n)_{n=1}^{\infty})=(y_n)_{n=1}^{\infty} \mbox{ if and only if } \hat{L}(S)((z_n)_{n=1}^{\infty})=(w_n)_{n=1}^{\infty} \] for every $(x_n)_{n=1}^{\infty}\in \mathcal{M}_{r}(\xi^1)$ where $(z_n)_{n=1}^{\infty}$ and $(w_n)_{n=1}^{\infty}$ are given by the above construction. The previous discussion implies that $\hat{L}(S)$ is defined at every element of $\mathcal{M}_{r}(L(\xi^1))$. 
	
	We now show that $\hat{L}(S)\in \mathcal{L}_T(L(\xi^1),L(\xi^2))$. Let $(z_n)_{n=1}^{\infty}\in \mathcal{M}_{r}(L(\xi^1))$ be arbitrary. It follows that \[ (\hat{T}_{L(\xi^2)}\circ \hat{L}(S))(z_n)_{n=1}^{\infty}=\hat{T}_{L(\xi^2)} (w_n)_{n=1}^{\infty}=(\xi_{n+1}^2Tw_n)_{n=1}^{\infty}\] where $(w_n)_{n=1}^{\infty}$ is such that $S(\xi_1^1z_1,z_1,z_2,\ldots)=(\xi_1^1w_1,w_1,w_2,\ldots)$. On the other hand, 
	\[
	(\hat{L}(S)\circ \hat{T}_{L(\xi^1)})(z_n)_{n=1}^{\infty}=\hat{L}(S)(\xi_{n+1}^1Tz_n)_{n=1}^{\infty}=(w_n')_{n=1}^{\infty} \] where $(w_n')_{n=1}^{\infty}$ is an abstract martingale satisfying	
	\begin{align*}
	(\xi_1^1 w_1', w_1',w_2',w_3',\ldots)&=S(\xi_1^1\xi_2^1 Tz_1,\xi_2^1 Tz_1,\xi_3^1 Tz_2,\ldots) \\ &=(S\circ \hat{T}_{\xi^1})(z_1,z_2,z_3,\ldots) \\ &=(\hat{T}_{\xi^2}\circ S)(\xi_1^1z_1,z_1,z_2,\ldots) \\ &=\hat{T}_{\xi^2} (\xi_1^1w_1,w_1,w_2,\ldots)\\ &=(\xi_1^2T\xi_1^1w_1,\xi_2^2Tw_1,\xi_3^2Tw_2,\ldots) \\ &=(\xi_1^1w_1',w_1',w_2',\ldots)
	\end{align*}
	Hence, $\hat{L}(S)\in \mathcal{L}_T(L(\xi^1),L(\xi^2))$.
\end{proof}

{\bf Acknowledgments} 
We would like to thank the referee for her/his valuable suggestions.

\end{document}